\documentclass[12pt]{amsart}
\usepackage[margin=0.6in]{geometry}

\usepackage[english]{babel}
\usepackage[alphabetic]{amsrefs}
\usepackage{amsmath,amssymb,amsfonts,amsthm,enumerate}
\usepackage{url}
\usepackage{graphicx}
\usepackage{wrapfig}

\usepackage[cal=boondox]{mathalfa}

\numberwithin{equation}{section}

\newtheorem{theorem}{Theorem}

\newtheorem{proposition}[theorem]{Proposition}

\newcommand{\R}{\mathbb{R}}
\newcommand{\Z}{\mathbb{Z}}
\newcommand{\N}{\mathbb{N}}

\renewcommand{\leq}{\leqslant}
\renewcommand{\le}{\leqslant}
\renewcommand{\geq}{\geqslant}
\renewcommand{\ge}{\geqslant}

\renewcommand{\epsilon}{\varepsilon}

\newcommand{\e}{\varepsilon}

\newcommand{\lam}{\lambda}

\begin{document}

\title{A fractional glance to the theory of edge dislocations}

\author{Serena Dipierro}
\address[Serena Dipierro]{Department of Mathematics and Statistics,
University of Western Australia,
35 Stirling Highway,
Crawley WA 6009, Australia}
\email{serena.dipierro@uwa.edu.au}

\author{Stefania Patrizi}
\address[Stefania Patrizi]{Department of Mathematics,
The University of Texas at Austin, 2515 Speedway Austin TX, 78712, USA} 
\email{spatrizi@math.utexas.edu}

\author{Enrico Valdinoci}
\address[Enrico Valdinoci]{Department of Mathematics and Statistics,
University of Western Australia,
35 Stirling Highway,
Crawley WA 6009, Australia}
\email{enrico.valdinoci@uwa.edu.au}

\begin{abstract}
We revisit some recents results inspired by the Peierls-Nabarro model on edge dislocations for crystals which rely on the fractional Laplace representation of the corresponding equation. In particular, we discuss results related to heteroclinic, homoclinic and multibump patterns for the atom dislocation function, the large space and time scale of the solutions of the parabolic problem, the dynamics of the dislocation points and the large time asymptotics after possible dislocation collisions.
\end{abstract}

\maketitle

\section{Introduction}

\subsection{The Peierls-Nabarro model for edge dislocations}\label{INTR}

The theory of edge dislocations describes defects
in the atoms display of a crystal structure.
In a nutshell, dislocations represent areas of a crystal
in which the atoms do not lie in the geometrically
organized pattern of the material,
and they are responsible
for plastic deformations to occur.
In particular, these plastic deformations
occur without altering
the chemical composition of the material and under external forces
significantly lower than expected.

\begin{figure}[h]
\centering
\includegraphics[height=3.8cm]{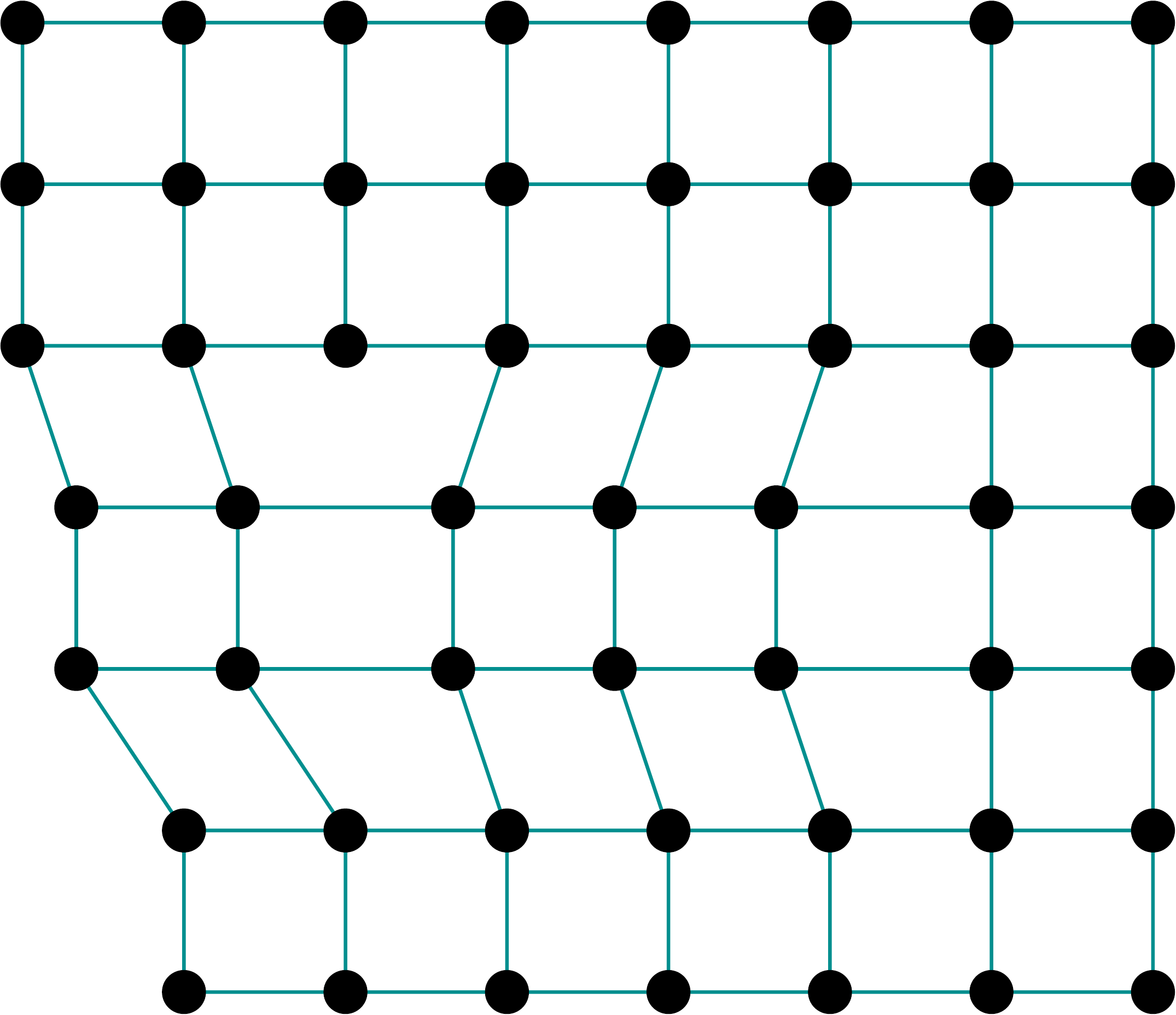}
\caption{A defect in a crystal (the black dots corresponding to
atoms and the cyan segments to the bonds between atoms).}\label{DEFE}
\end{figure}

As an example of dislocation, one can think of a simple cubic lattice representing
the crystal at a large scale and suppose that an irregularity arises
in the lattice: see for instance Figure~\ref{DEFE},
in which an imperfection of the crystal is visible on the fourth horizontal row.

As a matter of fact, Figure~\ref{DEFE} must be understood as (the local reproduction of) a planar section
of a three-dimensional (infinitely extended) crystal and represents the ``typical'' anomaly associated to an edge dislocation,
as being caused by the termination of a plane of atoms in the middle of a crystal
(the picture is similar to that of sticking half of a piece of paper into a stack of paper; in this case the
half of a piece of paper corresponding to the three extra atoms visible in the upper part of the third column of atoms in
Figure~\ref{DEFE}).

A pioneer description of elastic dislocations dates back to Vito Volterra~\cite{MR1509085}
and the connection between plastic deformations and dislocation theory
was set forth by Egon Orowan, Michael Polanyi and Geoffrey Taylor (independently and at about the same time)
explaining that the stress needed for the slip was far lower than theoretical predictions
in view of the theory of dislocations~\cite{ORO1, ORO2, ORO3, POLA, zbMATH02539109, zbMATH02539110}.
An intuitive reason for why dislocation dynamics allow plastic displacements under a low external force
can be found by comparing this phenomenon with 
the amusing locomotion of caterpillars
obtained by squeezing muscles in sequence in an undulating wave motion:
namely, since caterpillars do not have any bone in their entire bodies,
they can produce a ``defect'' by arching their body and obtain a displacement by moving this defect
along their body (this is however a rather crude biological explanation,
if interested in a detailed analysis of the caterpillars' locomotion see e.g.~\cite{VAN}).

At that time, the description of dislocations was mainly based on theoretical predictions
with little or no experimental evidence for dislocations. Indeed, although slip lines on the surface of metals had been
observed by optical microscopy,
a direct observation of crystal dislocations arrived later, since it
required transmission electron microscopes~\cite{VERMA, AME, BOL, HIR0, Silcox1959-SILDOO-2, HIR}.

\begin{figure}[h]
\centering
\includegraphics[height=3.8cm]{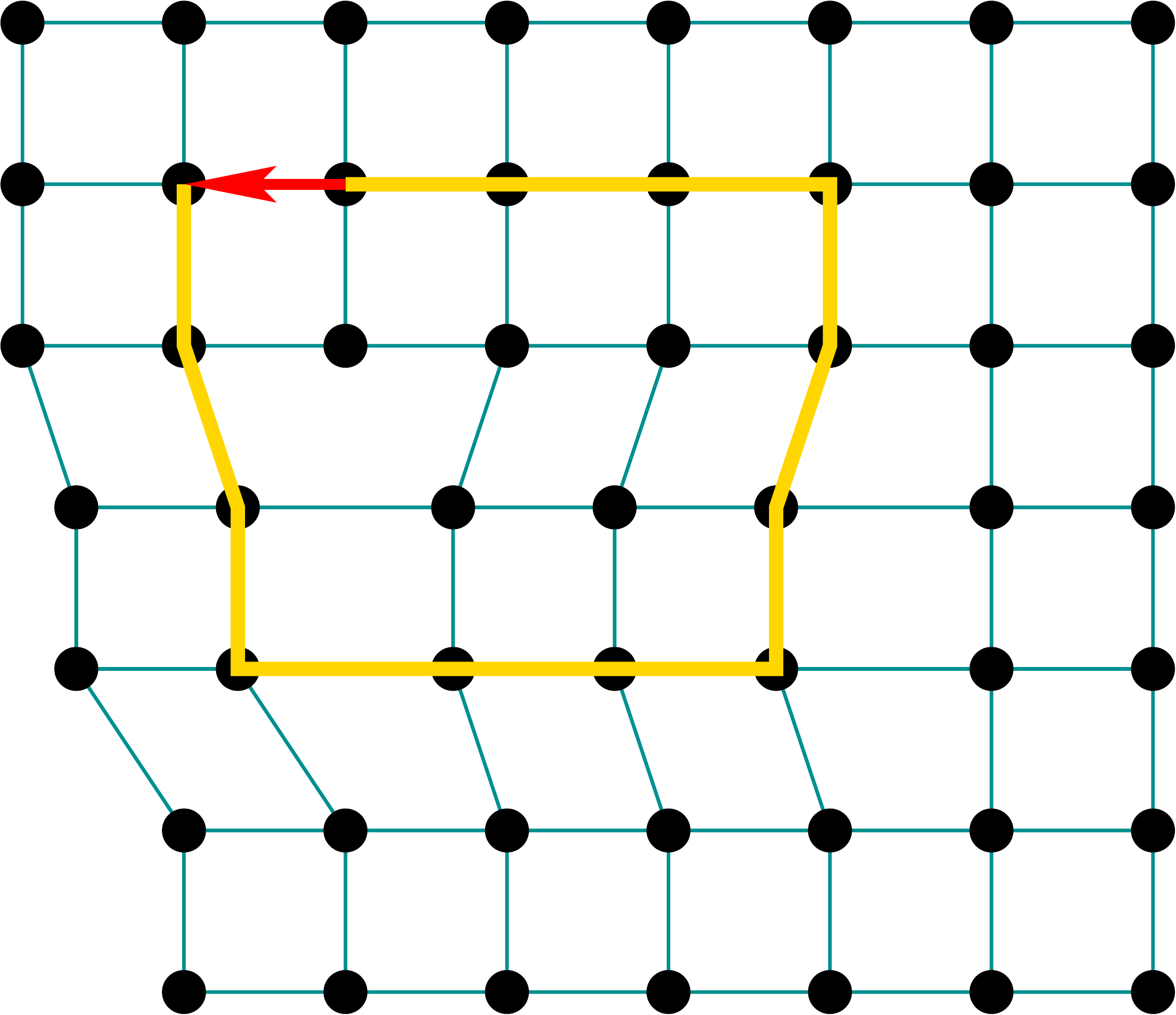}
\caption{Burgers vector.}\label{DEFF}
\end{figure}

A useful graphical device to detect crystal defects in a lattice is given by the so-called
Burgers vector, named after Johannes (Jan) Martinus Burgers~\cite{BUR}, see Figure~\ref{DEFF}.
The gist is that in a perfect crystal structure if one moves for a given number of steps, say three,
in every coordinate direction, say right-down-left-up, the trajectory produced is a square that
returns to the initial point; instead, when a similar ideal movement is
drawn encompassing a dislocation, the arrival point is different from the initial one
and the vector joining the starting position with the final one
is a vector, namely the Burgers vector, which encodes the magnitude and direction
of the crystal imperfection.

{F}rom the static viewpoint, a crystal defect produces a structural inhomogeneity
corresponding e.g. in Figure~\ref{DEFE} to a relative compression above the dislocation,
a relative tension below and a shear on the side (the shear can be visualized since
the bonds between atoms are shaped in the form of a parallelogram rather than a regular square).

\begin{figure}[h]
\centering
\includegraphics[height=3.8cm]{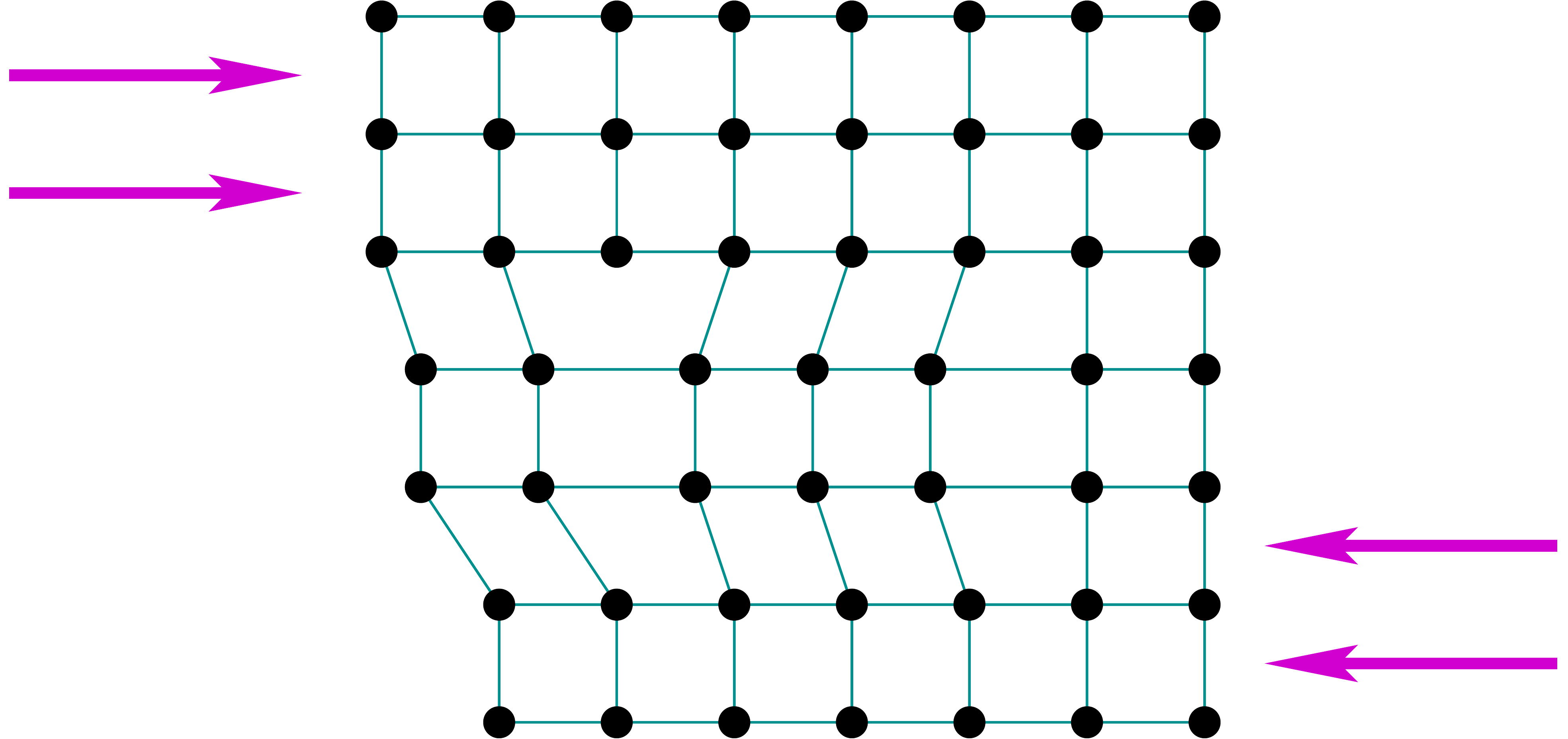}
\caption{Shearing forces.}\label{DEFG}
\end{figure}

In practice, one can take advantage of this structural imperfection
by producing plastic deformations as an outcome of the motion of dislocations
which is induced by a set of shearing forces,
namely forces parallel to one of the axis of the lattice,
pushing the crystal in one specific direction
from a side, and in the opposite direction from the other side
(whatever ``side'' may mean here, since our ideal crystal has infinite extension in every direction).
See Figure~\ref{DEFG} for a description of shearing forces.

\begin{figure}[h]
\centering
\includegraphics[width=0.2\textwidth]{fig1.pdf} $\mapsto$
\includegraphics[width=0.2\textwidth]{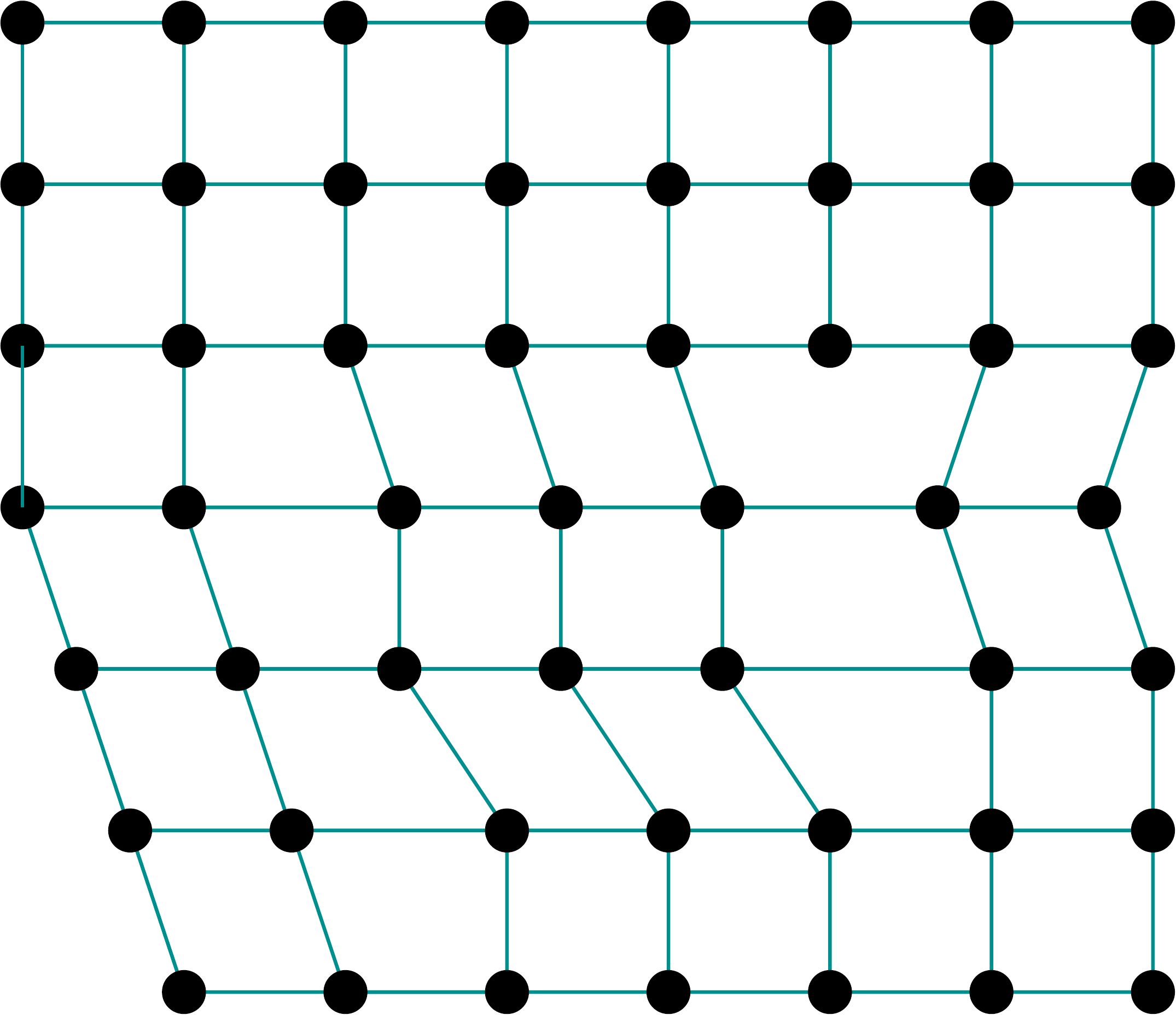} $\mapsto$
\includegraphics[width=0.2\textwidth]{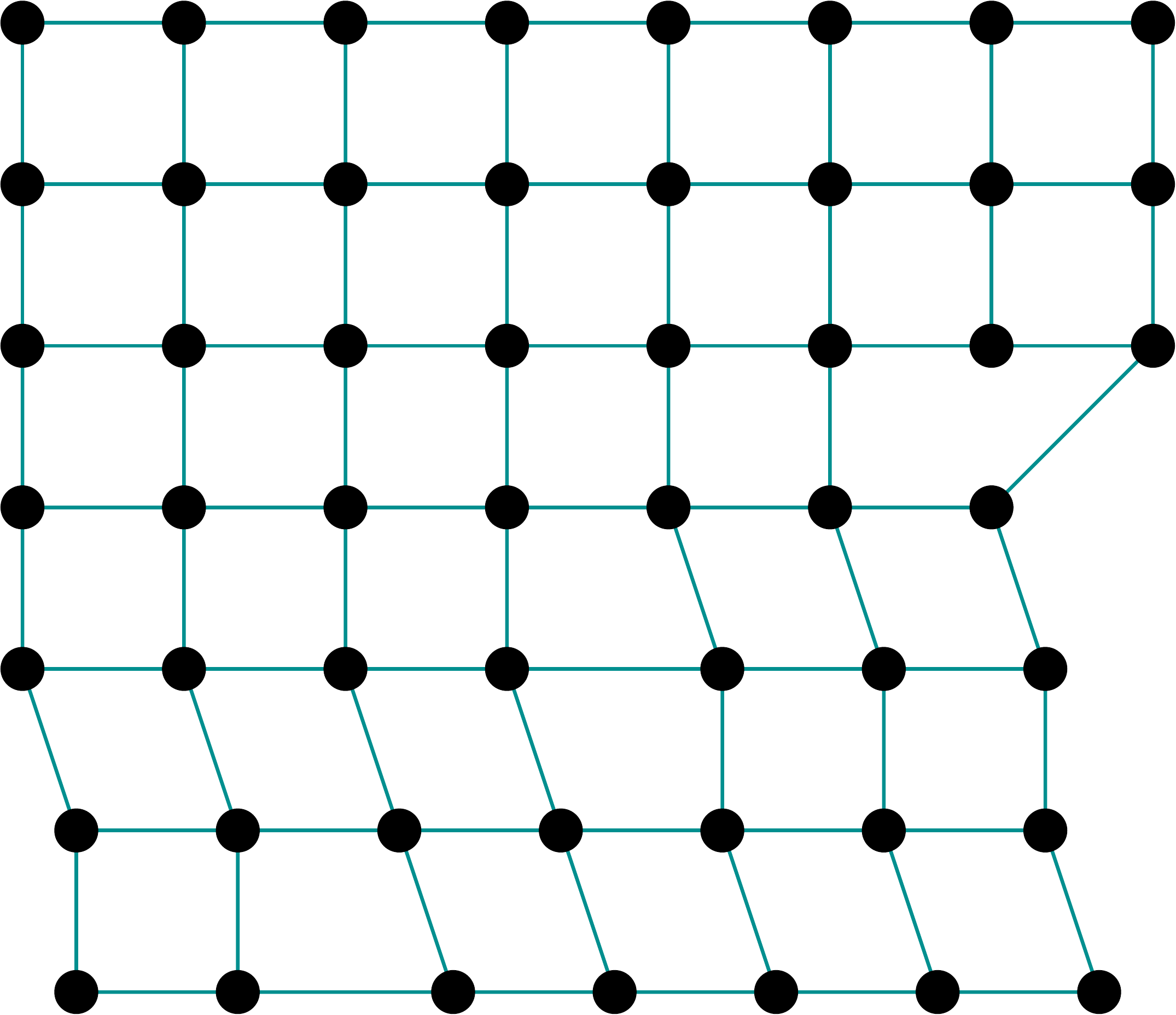} $\mapsto$
\includegraphics[width=0.2\textwidth]{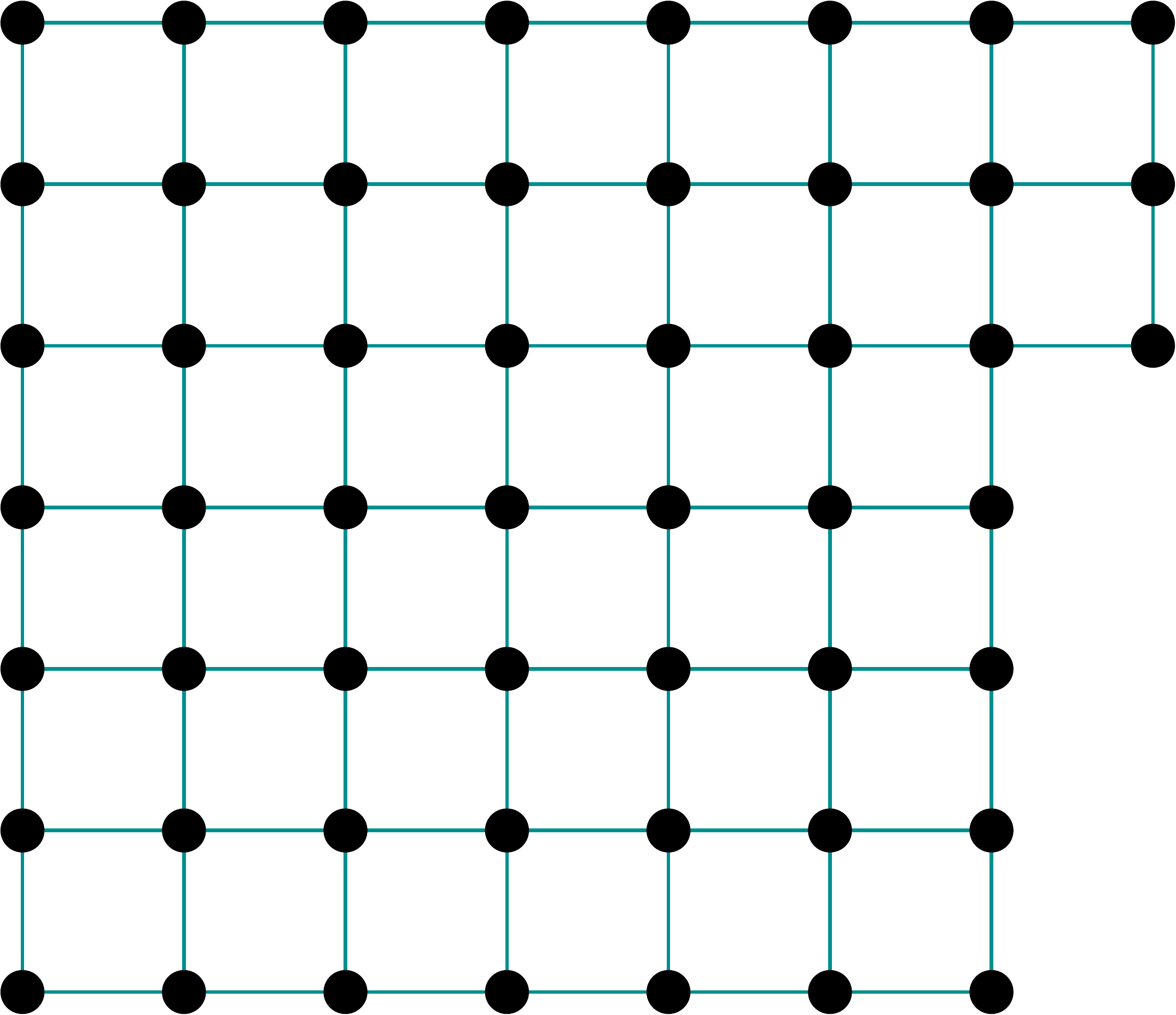} 
\caption{Plastic deformation as an outcome of crystal dislocation motion.}\label{DEFH}
\end{figure}

The advantage of dislocations towards plastic deformations is that when the
shearing forces are parallel to the Burgers vector their result can be that of weakening or breaking
some of the bonds between atoms and allow them to slide over each other.
This glide, or slip, method is thus capable of producing plastic deformations
at low stress levels. See Figure~\ref{DEFH} for a cartoon of
how a shearing force can move around dislocations in the lattice with the final outcome
of producing a visible plastic deformation:
the idea is that by conveniently sliding a dislocation
the crystalline order is restored on either side but the atoms on one side have moved by one position.

\begin{figure}[h]
\centering
\includegraphics[height=4.5cm]{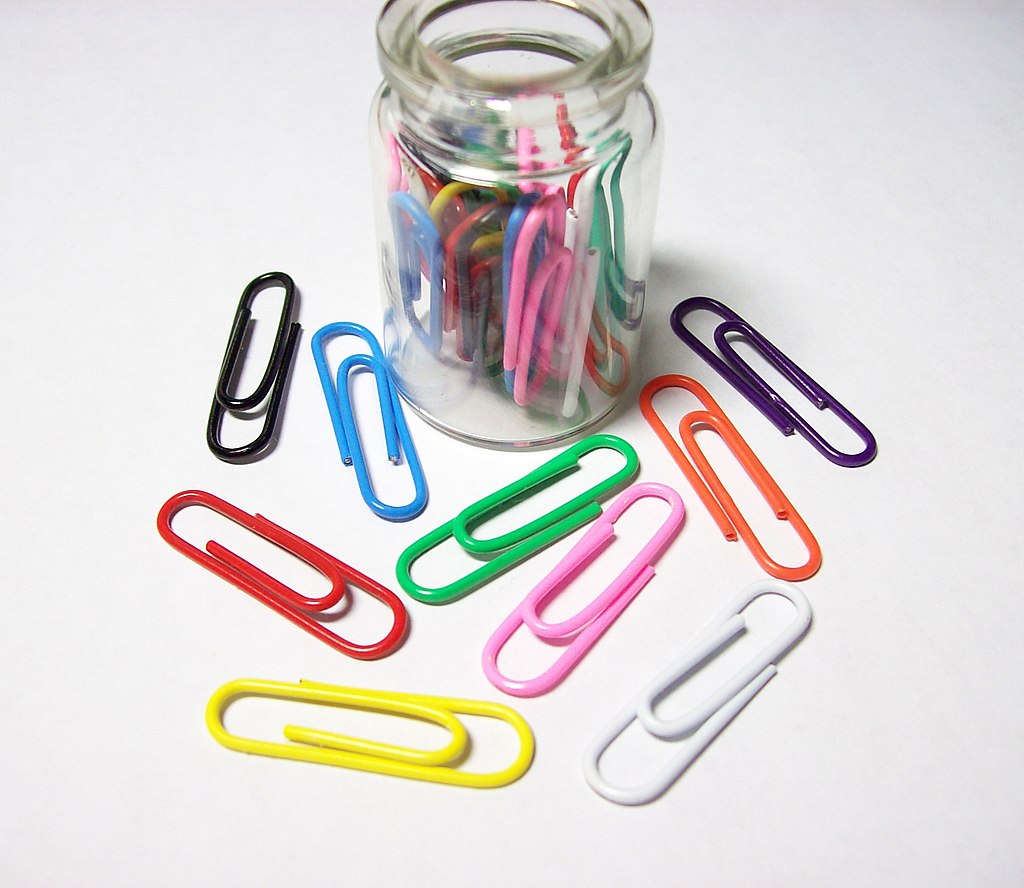}
\caption{A few paper clips of different colors (Public Domain image from Wikipedia).}\label{DEFG8}
\end{figure}

We have a direct experience of the role played by dislocations
in plastic deformations anytime we indulge in playing around with paperclips
(those little devices made of steel wire bent to a looped shape
that we use to hold sheets of paper together before every useful information was digitalized
and posted on the internet, see Figure~\ref{DEFG8}): with soft movements of our fingers,
we can easily open and close repeatedly the looped shape of a paperclip,
enjoying its flexibility and malleability -- till at some point the paperclip becomes suddenly quite
stiff and to change its shape we need to exert some extra force which typically ends
up breaking the paperclip. Hence, this simple experiment somehow confirms
how the strength and ductility of metals are influenced and controlled by
dislocations.

The physical explanation of this phenomenon is that at the beginning
the paperclip structure presents some defects and our soft movements
rely on the dynamics of these dislocations to modify the shape of the object;
at some point, however, several of these imperfections get resolved precisely by our
repeated opening and closing the loop of the clip, and when no more dislocation is available
(or when too few dislocations are available) the perfect (or almost perfect) lattice
offers a sounder resistance to the external stress (see the last frame in Figure~\ref{DEFH}), which results in the (perhaps not so expected)
stiffness that typically make this game end by the breaking of the clip.

\begin{figure}[h]
\centering
\includegraphics[height=5cm]{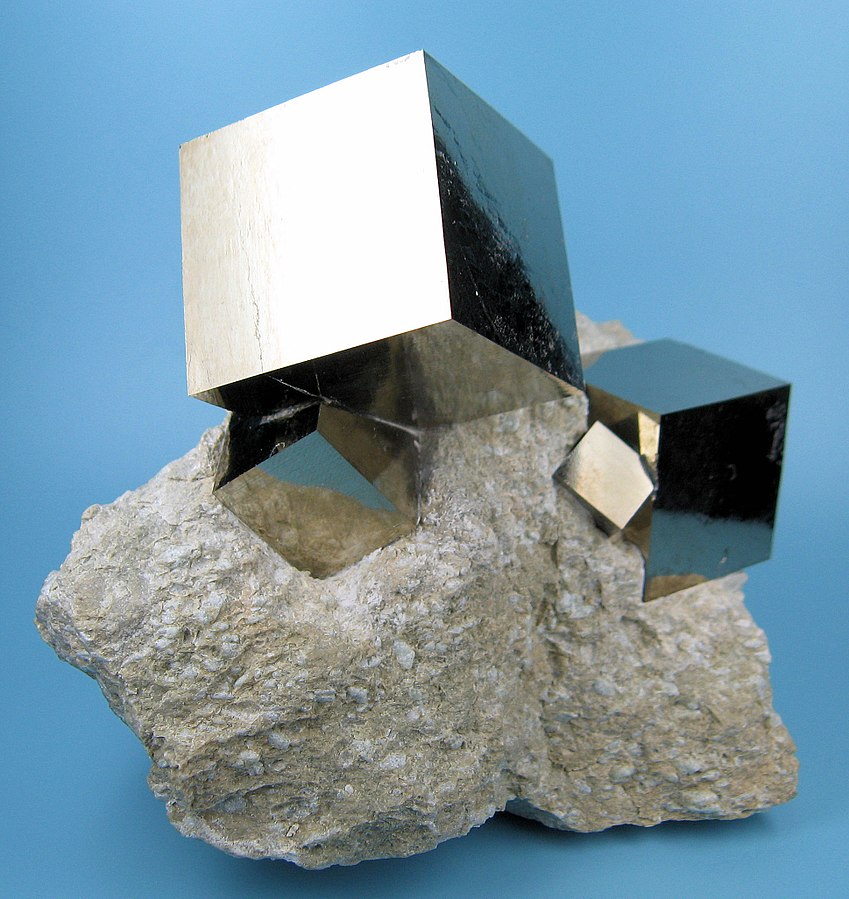}
\caption{Pyrite cubic crystals (image from Wikipedia,
licensed under the Creative Commons Attribution-Share Alike 3.0 Unported license).}\label{DEFG9}
\end{figure}

Some technical remarks are in order. First off, the assumption that the crystal has a simple cubic structure
(also called primitive cubic in jargon) is indeed a great geometric simplification
but it includes several concrete cases of interest, such as 
polonium~${\rm Po}$, see~\cite{PLON}, iron pyrite~${\rm FeS_2}$,
see Figure~\ref{DEFG9}, as well as an allotrope  of selenium~${\rm Se}$,
see~\cite{SELE}.

\begin{figure}[h]
\centering
\includegraphics[width=0.22\textwidth]{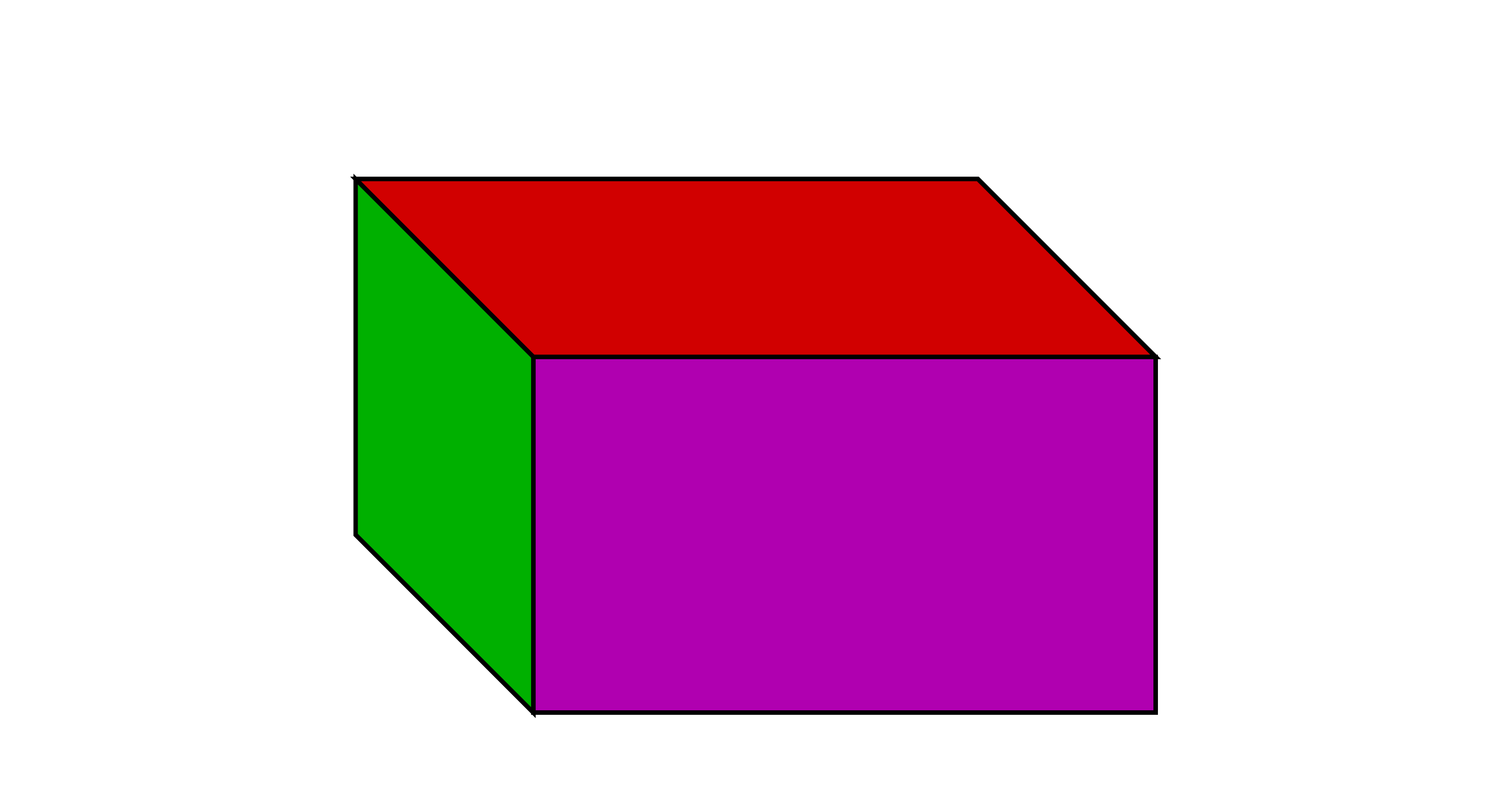} $\mapsto$
\includegraphics[width=0.22\textwidth]{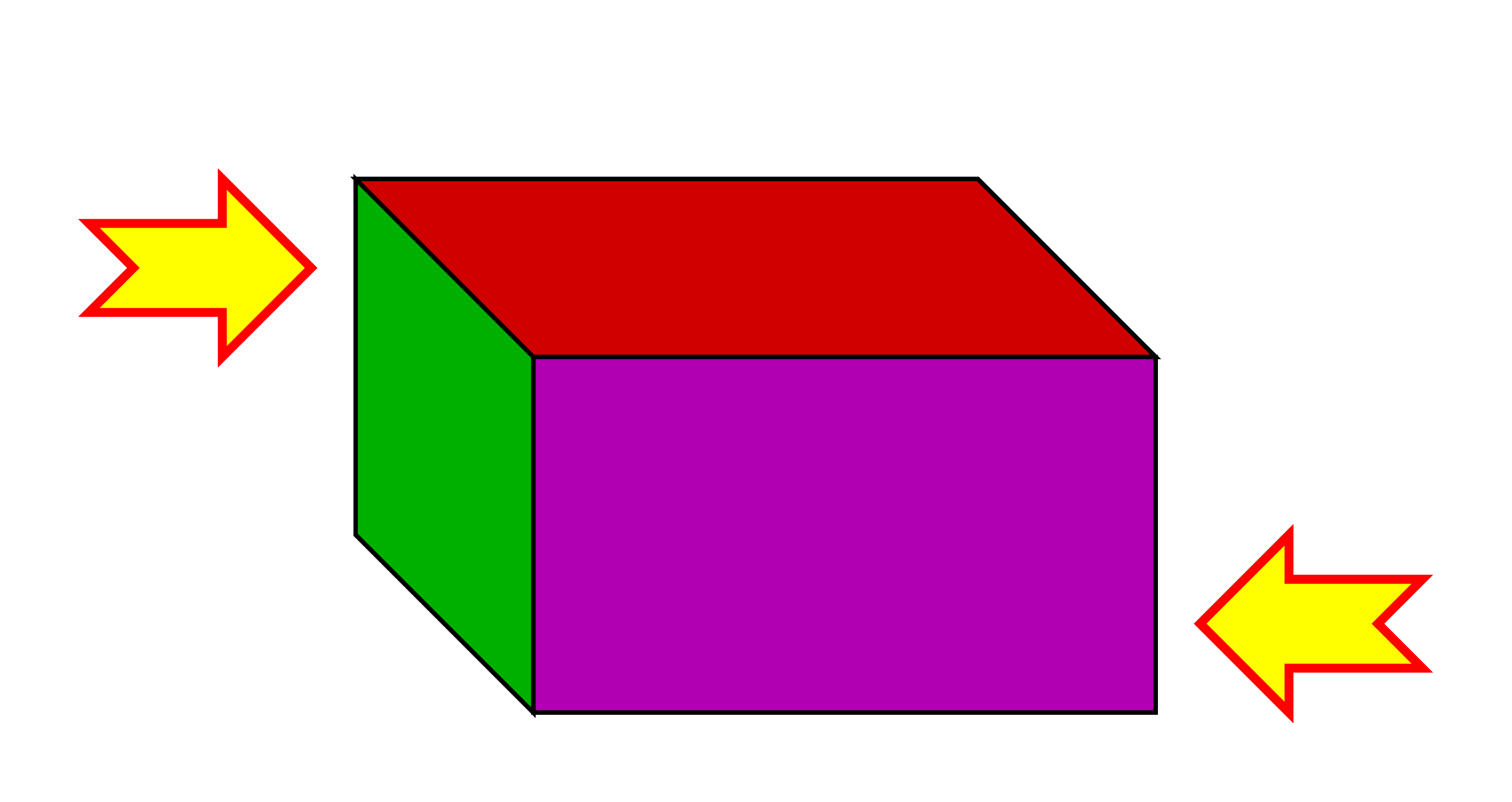} $\mapsto$
\includegraphics[width=0.22\textwidth]{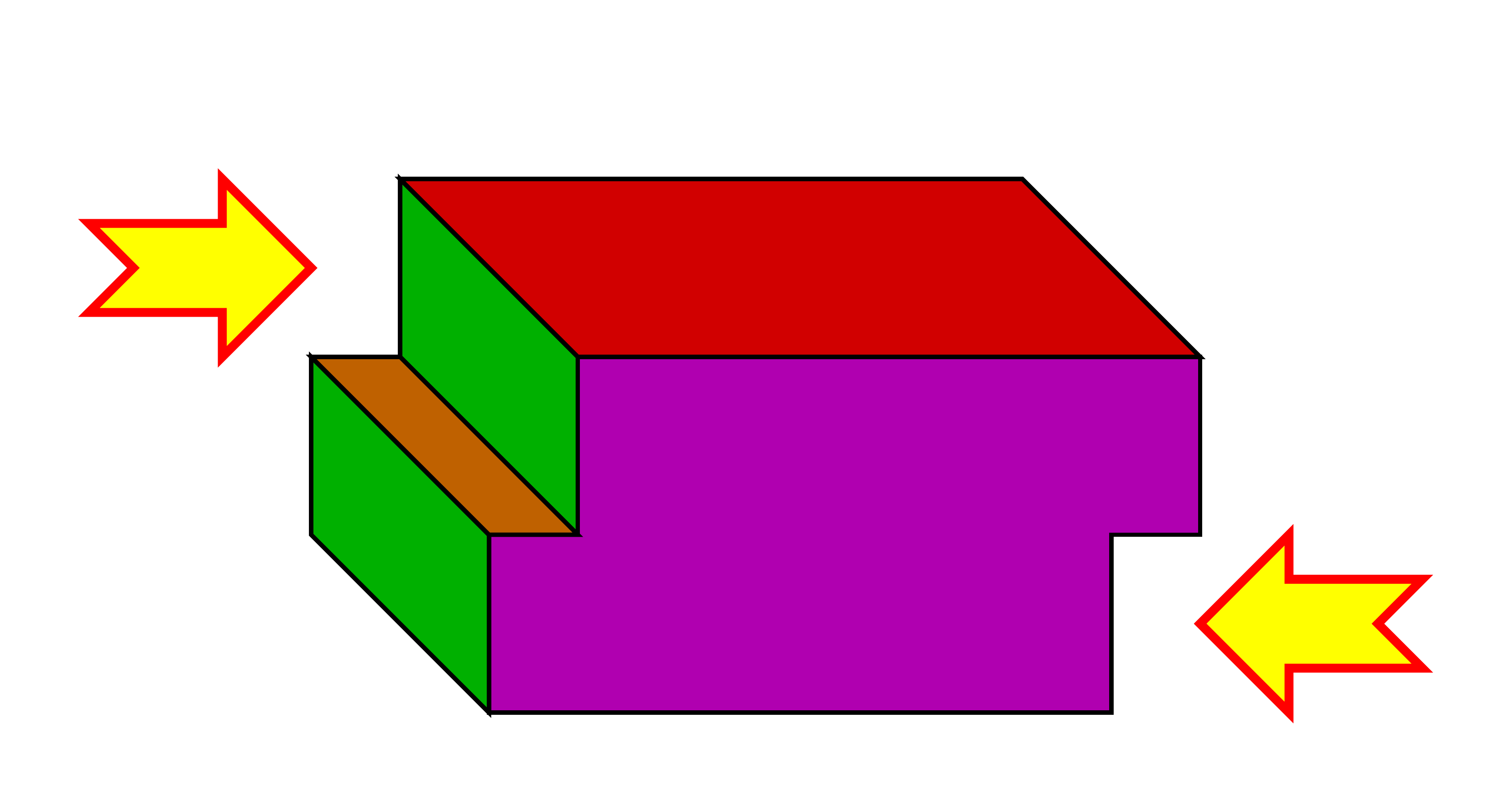} $\mapsto$
\includegraphics[width=0.22\textwidth]{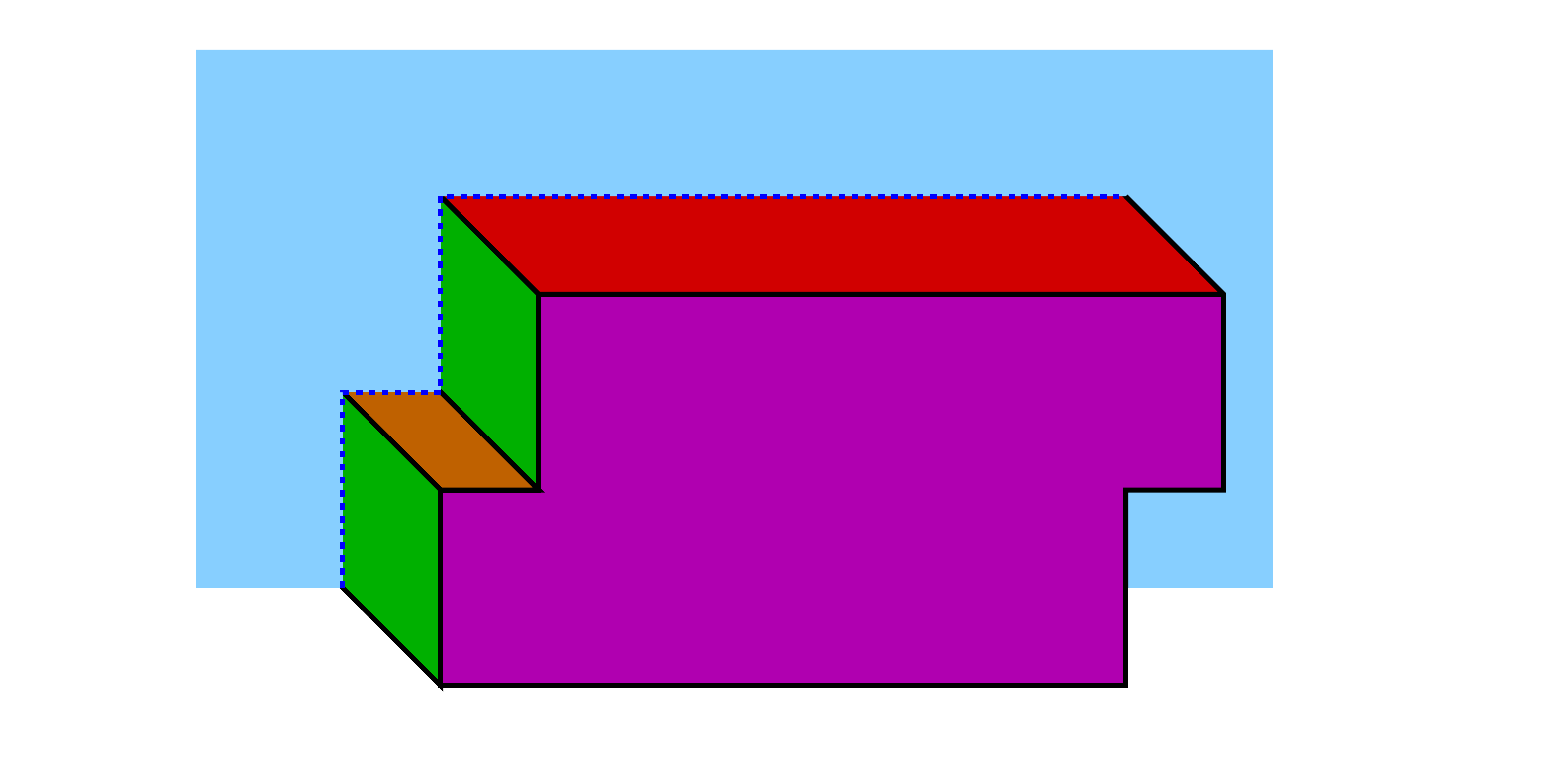} 
\caption{Slip of a crystal along a horizontal plane in an edge dislocation, and vertical section.}\label{DEFHSe}
\end{figure}

\begin{figure}[h]
\centering
\includegraphics[width=0.22\textwidth]{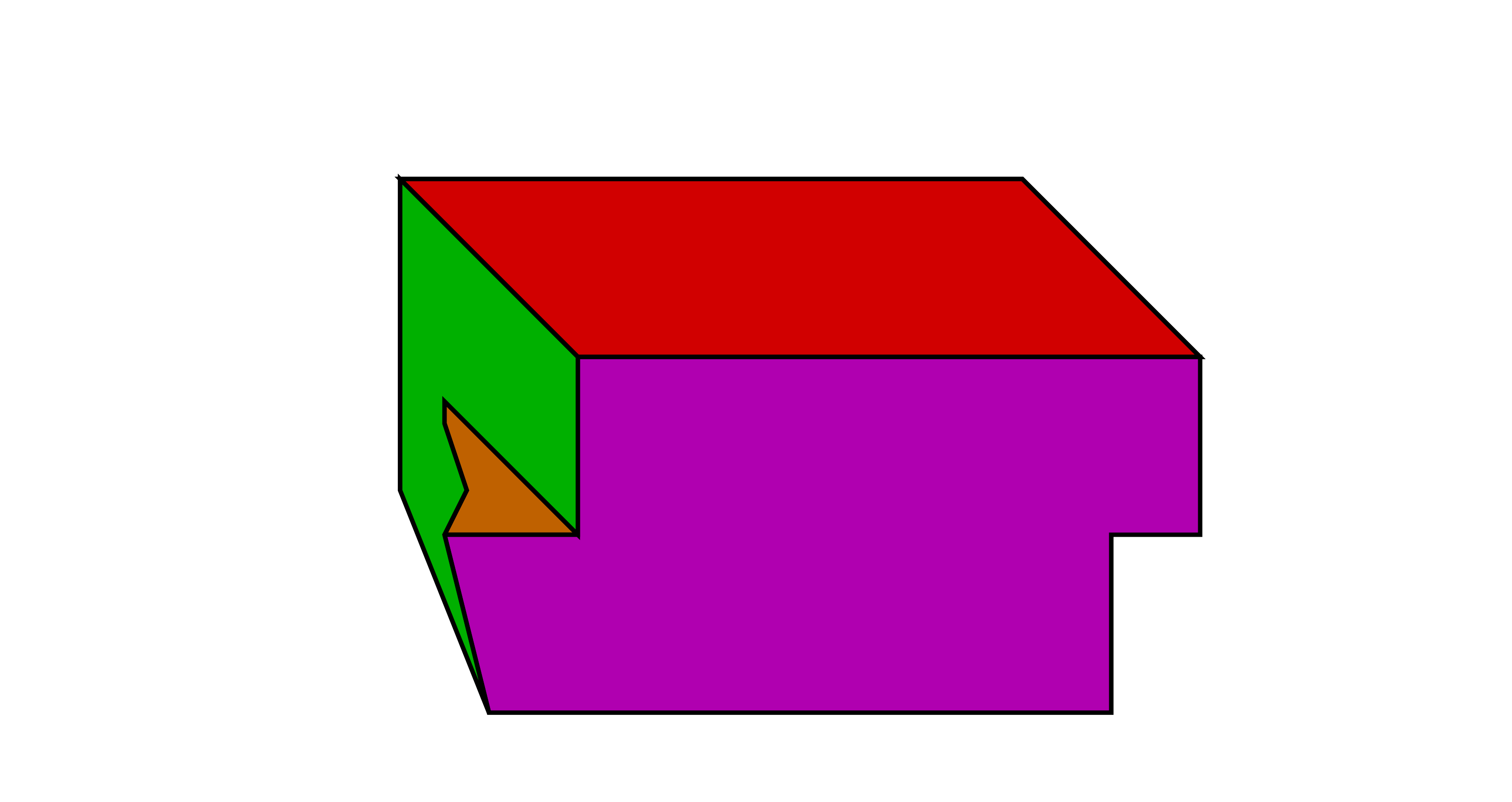}
\caption{A screw dislocation.}\label{DEFHSe45}
\end{figure}

Furthermore, we stress that we are focused here on edge dislocations,
namely we only consider here situations
in which dislocations move in the direction of the Burgers vector
(as in Figures~\ref{DEFG}, \ref{DEFH} and~\ref{DEFHSe}).
Other situations also occur in nature, such as the so-called
screw dislocations in which the dislocation dynamics takes place in a direction perpendicular to the Burgers vector, thus creating a helical arrangement of atoms around the core,
see e.g. Figure~\ref{DEFHSe45}.
Mixed dislocations exist as well, in which case the Burgers vector is neither parallel
nor orthogonal to the dislocation dynamics. But in this note we confine our interest to the simplest
case of edge dislocations.

To mathematically describe edge dislocations, we recall the classical model 
by Rudolf Peierls and Frank Nabarro. Suppose that the crystal is subject to shearing forces
causing a slip along a horizontal plane, that we call the slip plane. We take a vertical section of the crystal,
see Figure~\ref{DEFHSe}.
The intersection between the horizontal slip plane and the vertical sectional plane is called slip line.
We suppose that at a large scale the crystal has a simple cubic structure described
by the lattice~$\Z^3$ and accordingly the vertical sectional plane
is such that it detects a large scale atom disposition described by the lattice~$\Z^2$.

We stress that this lattice does not describe the position of every single atom,
since defects at small scales can occur. The dislocation theory thus makes the ansatz that
the lattice structure simply dictates a favorable state for atoms which are however free to move,
subject to a potential that penalizes the atom positions outside the lattice structure
and a bond interaction between atoms.\footnote{To favor the intuition, one can think about
people in a room with chairs placed in a lattice structure. Suppose that people has some bonds
between them (either psychological, being influenced by the movements of the people around
and trying to follow their neighbors, or physical, e.g. holding hands). Suppose also that
people have a preference to sit comfortably on a chair (not necessarily ``their'' chair, any chair
is equally good). Most chairs are occupied, but some chair is free; most people sit on a chair, but perhaps someone
could not find a convenient place.
Then, if some people moves around, they may drag with them some other people,
but the collective interest is to find a chair to sit down.

With respect to the crystal dislocation model,
people correspond to the atoms, the chair to the atom favorable position aligned with the cubic
lattice, the empty chairs to the crystal defects. A quantitative difference is that typically the
potential describing people's satisfaction is discontinuous (one would charge zero discomfort
if someone sits on a chair and unit discomfort if not, no matter how far or how close the individual is to
the closest chair), while the potential used in physical applications are typically of multi-well shapes
(the minima corresponding to the optimal location on the lattice), but continuous and in fact smooth.}

The slip will indeed cause some misalignment between the atoms below the slip lines and the ones above. For concreteness, we focus our attention on the atoms above the slip line
(the case of the ones below being similar) and we aim at describing their mismatch with respect
to the lattice structure.

\begin{figure}[h]
\centering
\includegraphics[width=0.77\textwidth]{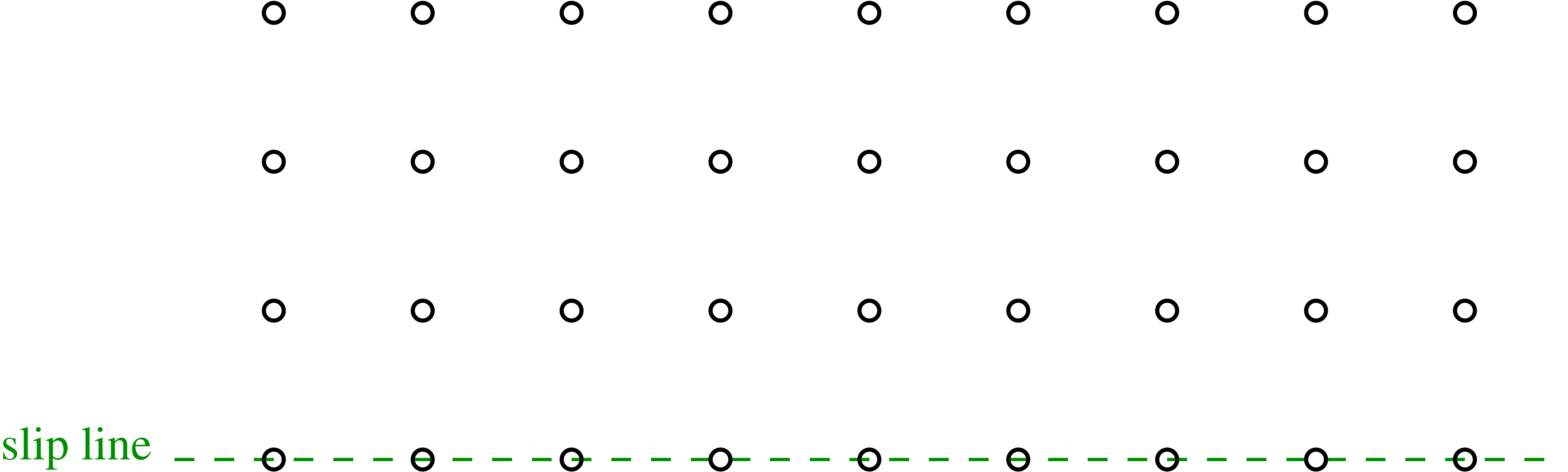}
\caption{Rest position of the atoms corresponding to the large scale crystalline structure.}\label{DEFHSfg1}
\end{figure}

For this, we consider the atom rest positions~$\Z^2$ in the upper halfplane~$\R\times[0,+\infty)$,
corresponding to the empty circles in Figure~\ref{DEFHSfg1}.

\begin{figure}[h]
\centering
\includegraphics[width=0.77\textwidth]{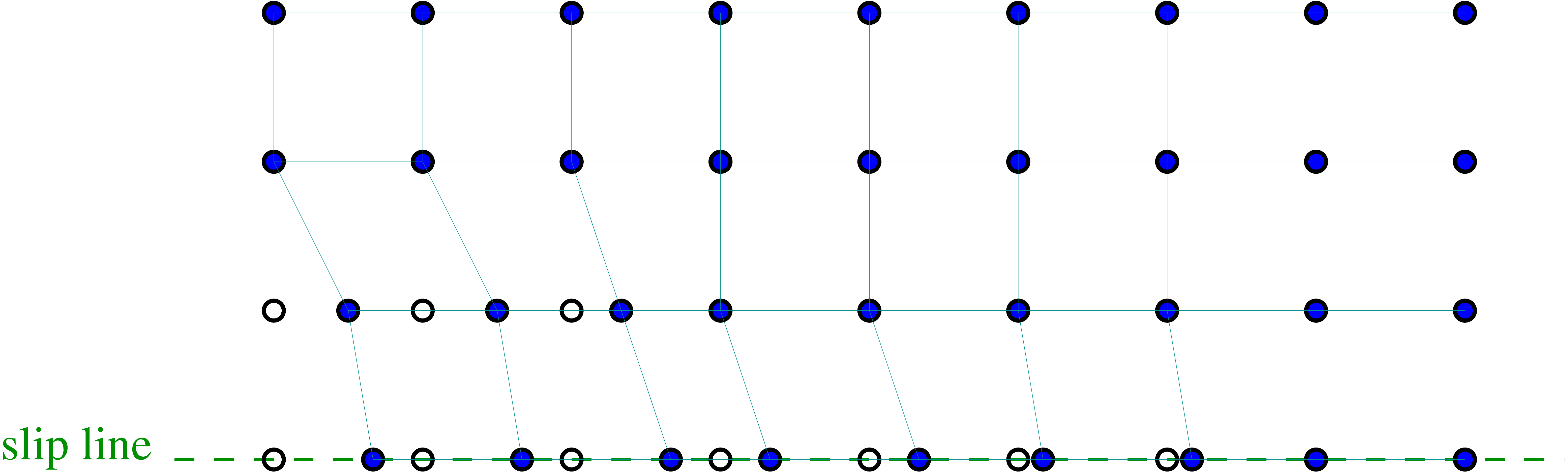}
\caption{Mismatch between atom location and crystal structure.}\label{DEFHSfg2}
\end{figure}

Positions in this halfplane are denoted by using coordinates~$(x,y)\in\R\times[0,+\infty)$.
The position of each atom is described by a dislocation function that
identifies the mismatch with respect to the rest position: for simplicity,
we assume that atoms can only slide horizontally (that is,
parallel to the slip plane and to the shearing forces), see Figure~\ref{DEFHSfg2}
in which the full circles represent the atoms (we are thus neglecting
the vertical mismatch of atoms possibly caused by mutual bonds, by considering
this effect as negligible with respect to the shearing force).

\begin{figure}[h]
\centering
\includegraphics[width=0.29\textwidth]{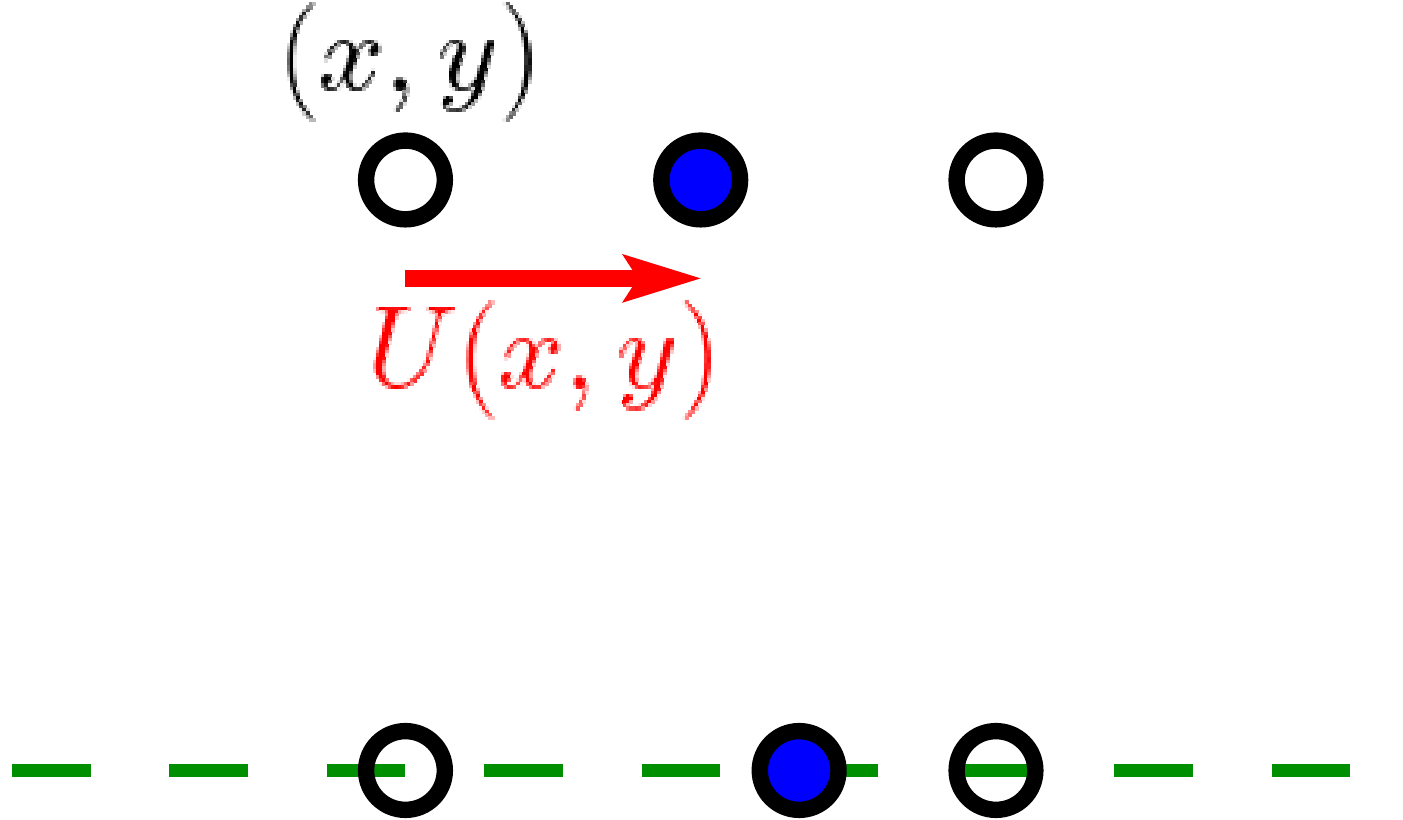} 
\caption{The dislocation function~$U(x,y)$.}\label{DEFHSfh}
\end{figure}

In this setting, the distance between the actual position of an atom
and its rest position  located at~$(x,y)\in\R\times[0,+\infty)$
corresponds to a scalar function~$U(x,y)$ (say, with~$U$
positive if the atom is dislocated towards the right \label{NOTDISDS}
and negative if the atom is dislocated towards the left, see Figure~\ref{DEFHSfh}).
It is also convenient to consider the trace of the dislocation function~$U$ along the slip line
which corresponds to~$u(x):=U(x,0)$.

In a situation such as the one described in Figure~\ref{DEFHSe}, we can argue that
most of the mismatch occurs along the slip line: roughly speaking, we can expect
the atoms located in the upper halfplane to be mostly alligned with a cubic lattice
and the ones located in the lower halfplane to be mostly alligned with another cubic lattice
which is a horizontal translation of the previous one.
Therefore, we can suppose that the main contribution to the crystal mismatch
is encoded by the function~$u$. To quantify this mismatch, we consider a multiwell potential~$W$,
say a scalar function such that~$W(k)=0$ for all~$k\in\Z$ and $W(r)>0$ for all~$r\in\R\setminus\Z$.
For simplicity, we will suppose in what follows that~$W$ is smooth,
periodic of period~$1$, even and such that~$W''(0)\ne0$.

This potential penalizes atoms which are not aligned with the lattice structure, since~$W(u(x))$
turns out to be positive unless~$u(x)\in\Z$, that is unless the atom with ideal rest position at~$x$
is located to a point of the lattice (recall that any point of the lattice ends up to be good
for the atom to be in equilibrium with respect to the large scale structure of the crystal).

As a result, in this model the total energy contribution coming from atom displacement can be written as
$$ {\mathcal{P}}(U):=\int_\R W(u(x))\,dx.$$
The total energy however has to account also for the bonds between atoms.
The simplest ansatz is to suppose that this energy contribution resembles an elastic energy of the form
\begin{equation}\label{SQRL4} {\mathcal{E}}(U):=\frac12\int_{\R\times(0,+\infty)} |\nabla U(x,y)|^2\,dx\,dy,\end{equation}
where the factor~$\frac12$ has been introduced for later convenience (we are disregarding here
structural constants).

The total energy is therefore the sum of the displacement energy and the bond energy and takes the form
$$ {\mathcal{H}}(U):={\mathcal{E}}(U)+
{\mathcal{P}}(U)=\frac12\int_{\R\times(0,+\infty)} |\nabla U(x,y)|^2\,dx\,dy
+\int_\R W(u(x))\,dx.$$
We stress that this model is, in a sense, of ``hybrid'' type:
while the rest position of the atoms correspond to a lattice structure with a discrete feature,
the dislocation function~$U$ (and thus its trace~$u$)
is assumed to be continuous (and actually smooth) and defined
in the continuum given by~$\R\times(0,+\infty)$.
This ``philosophical inconsistency'' between the discrete and continuous descriptions of nature
is compensated by a simple and effective realization of the model
and it is also, at least partially, justified by the fact that
the atom distance is typically much smaller than the plastic effects that one is interested in taking into
account (therefore, for practical purposes, the lattice acting as a substratum is
``almost'' a continuum).

Equilibrium configurations of this model correspond to critical points
of~${\mathcal{H}}$.
These configurations can be detected by noticing that, for every~$R>0$
and any smooth function~$\Psi:
\R\times[0,+\infty)\to\R$ such that~$\Psi(x,y)=0$ if~$|(x,y)|\ge R$,
\begin{equation}\label{HHJA}
\begin{split} 0&=
\langle D{\mathcal{H}}(U),\Psi\rangle\\&=
\int_{\R\times(0,+\infty)} \nabla U(x,y)\cdot\nabla\Psi(x,y)\,dx\,dy
+\int_\R W'(u(x))\psi(x)\,dx\\
&=-\int_{\R\times(0,+\infty)} \Delta U(x,y)\Psi(x,y)\,dx\,dy
+\int_\R \big(W'(u(x))-\partial_y U(x,0)\big)\psi(x)\,dx,
\end{split}\end{equation}
where~$\psi(x):=\Psi(x,0)$, thanks to Green's first identity.

In particular, if we take a ball~$B\Subset\R\times(0,+\infty)$ and a function~$\Psi_0\in C^\infty_0(B)$,
it follows that~$\psi_0(x):=\Psi_0(x,0)=0$ for all~$x\in\R$ and thus, by~\eqref{HHJA},
$$ 0=-\int_{\R\times(0,+\infty)} \Delta U(x,y)\Psi(x,y)\,dx\,dy,$$
that leads to~$U$ being harmonic in~$\R\times(0,+\infty)$.

Inserting this information into~\eqref{HHJA}, we find that
\begin{equation*}
0=\int_\R \big(W'(u(x))-\partial_y U(x,0)\big)\psi(x)\,dx,
\end{equation*}
giving that~$W'(u)=\partial_y U$ along the slip line.

{F}rom these considerations, we infer that equilibria of the 
Peierls-Nabarro model correspond to solutions of
\begin{equation}\label{SQRL2}
\begin{cases}
\Delta U=0 & {\mbox{ in }}\R\times(0,+\infty),\\
\partial_y U=W'(u) & {\mbox{ in }}\R\times\{0\}.
\end{cases}
\end{equation}
Formally, using the harmonicity of~$U$ and the Poisson kernel for the halfspace, we thus obtain that,
for all~$(x,y)\in\R\times(0,+\infty)$,
$$ U(x,y)=\frac1{\pi}\int_\R \frac{y\,u(\zeta)}{(x-\zeta)^2+y^2}\,d\zeta.$$
To remove the singularity as~$y\searrow0$, one can notice that a primitive
of~$\R\ni\zeta\mapsto \frac{y}{(x-\zeta)^2+y^2}$ is~$ \arctan\left(\frac{\zeta-x}y\right)$ and accordingly
$$ \frac1{\pi}\int_\R \frac{y}{(x-\zeta)^2+y^2}\,d\zeta
=\frac1{\pi}\big(\arctan(+\infty)-\arctan(-\infty)
\big)=1.
$$
This gives that
$$ U(x,y)-u(x)=\frac1{\pi}\int_\R \frac{y\,\big(u(\zeta)-u(x)\big)}{(x-\zeta)^2+y^2}\,d\zeta$$
and therefore
\begin{eqnarray*}&&W'(u(x))=\lim_{y\searrow0}
\partial_y U(x,y)=\lim_{y\searrow0}
\partial_y \big(U(x,y)-u(x)\big)=
\frac1{\pi}\lim_{y\searrow0}
\int_\R \frac{\partial}{\partial y}\left(\frac{y\,\big(u(\zeta)-u(x)\big)}{(x-\zeta)^2+y^2}\right)\,d\zeta\\&&\quad
=\frac1{\pi}\lim_{y\searrow0}
\int_\R\frac{ (\zeta^2 - 2 x\zeta + x^2 - y^2) \,\big(u(\zeta)-u(x)\big)}{( \zeta^2 - 2 x\zeta + x^2 + y^2)^2}\,d\zeta
=\frac1{\pi}
\int_\R\frac{ u(\zeta)-u(x)}{\zeta^2 - 2 x\zeta + x^2}\,d\zeta\\&&\quad=\frac1{\pi}
\int_\R\frac{ u(\zeta)-u(x)}{(x-\zeta)^2}\,d\zeta,
\end{eqnarray*}
that is, neglecting normalizing constants,
\begin{equation}\label{SQRL1}
-\sqrt{-\Delta }\,u=W'(u)\quad{\mbox{ in }}\,\R.
\end{equation}
Alternatively, one can formally obtain~\eqref{SQRL1} from~\eqref{SQRL2}
from the analysis of a Dirichlet-to-Neumann
problem (see e.g.~\cite[Section~4.1]{MR3469920})
or solving the linear equation in~\eqref{SQRL2} via Fourier transform
(see e.g.~\cite[Section~4.3]{MR3469920}).

We think that~\eqref{SQRL1} is a neat example of a nonlocal equation of fractional type
arising from purely classical considerations.
Of course, equation~\eqref{SQRL1} presents both advantages and disadvantages with respect to~\eqref{SQRL2}. On the one hand, equation~\eqref{SQRL2} is classical in nature
and it may seem more handy to be dealt with. On the other hand, equation~\eqref{SQRL1} showcases
the intrinsic nonlocal aspect of the problem, in which dragging one single atom along the slip line
ends up influencing, in principle, the whole crystalline structure. Furthermore, equation~\eqref{SQRL1}
is set in dimension~$1$, which is a considerable conceptual advantage
since it opens the possibility of employing methods from dynamical systems
(which will provide essential ingredients in Section~\ref{CHA}).

Additionally, equation~\eqref{SQRL1} immediately allows to consider more general equations of the form
\begin{equation}\label{SQRL3}
-(-\Delta )^s u=W'(u)\quad{\mbox{ in }}\,\R,
\end{equation}
for a given~$s\in(0,1)$. 

In this spirit, equation~\eqref{SQRL3} is not a mere generalization of~\eqref{SQRL1} motivated by mathematical curiosity: instead, it can offer a useful approach to comprise models
in which the range-effect of the atomic bonds is different from the one modeled by the purely elastic
response introduced in~\eqref{SQRL4}. Moreover,
it allows a direct comparison between equations arising in material sciences and crystal dislocation
dynamics with similar ones motivated, for instance, by problems in ethology, population dynamics and mathematical biology (for instance, similar models can be considered to describe
a biological population subject to a nonlocal dispersal strategy
of L\'evy flight type, see e.g.~\cite{GG}
and the references therein).

Also, and possibly more importantly, considering~\eqref{SQRL1}
as a ``particular case'' of~\eqref{SQRL3} allows one to take into account a continuous range of fractional
exponents~$s$, avoiding ``discrete jumps'' from the classical case~$s=1$
(often related to the classical Allen-Cahn equation) and the case~$s=\frac12$
in~\eqref{SQRL3}: besides its conceptual advantage, this approach allows one to exploit
perturbative methods in the fractional exponent (e.g. bifurcating from previously known cases).

\subsection{Organization of the paper}

The rest of the article is organized as follows.
In Section~\ref{CHA} we present some recent results that describe situations in which
atom dislocations in crystals can display complicated
or unusual patterns as stationary equilibria, thus producing chaotic behaviors.

In Section~\ref{DYN:SYS} we analyze the dislocation dynamics and in Section~\ref{COLL:SE}
the cases in which collisions between dislocation points occur.
Section~\ref{sec:rel0686uyf}
is then devoted to discussing what happens to the dislocations after collisions.

In Section~\ref{orowan123} we deduce the so-called 
Orowan's Law, and in Section~\ref{homog123} we deal with the homogenization problem
associated with the Peierls-Nabarro model.

\section{Chaotic dislocation displacements}\label{CHA}

A natural question is whether atom dislocations in crystals can display complicated
or unusual patterns as stationary equilibria. We started the investigation of this question
in~\cite{MR3594365, MR4053239}.
The gist of this analysis is that, without indulging in philosophical definitions,
one of the simplest forms of ``chaos'' consists in the possibility of producing patterns
which exhibit transition layers and oscillatory behaviors. To accomplish this goal,
we find it convenient to revisit some classical methods from dynamical systems and topological analysis
(such as the ones in~\cite{MR1119200, MR1799055})
in light of the fractional setting introduced in~\eqref{SQRL3}.

{F}rom the technical point of view, we recall that the existing methods
were designed for classical Hamiltonian or Lagrangian systems
and were aiming at detecting oscillatory behaviors in the time variable, that is
for the time evolution of a particle, or a system of particles. The structure in~\eqref{SQRL3}
is different, since no time variable is involved: instead, one can use methods from dynamical
systems to detect oscillatory phenomena in the space variable, especially thanks to the
fact that equation~\eqref{SQRL3} is set in one space variable (after all, once things
are set into a good mathematical formulation, one space variable can formally play the role of time,
and vice versa).

\begin{figure}[h]
\centering
\includegraphics[width=0.4\textwidth]{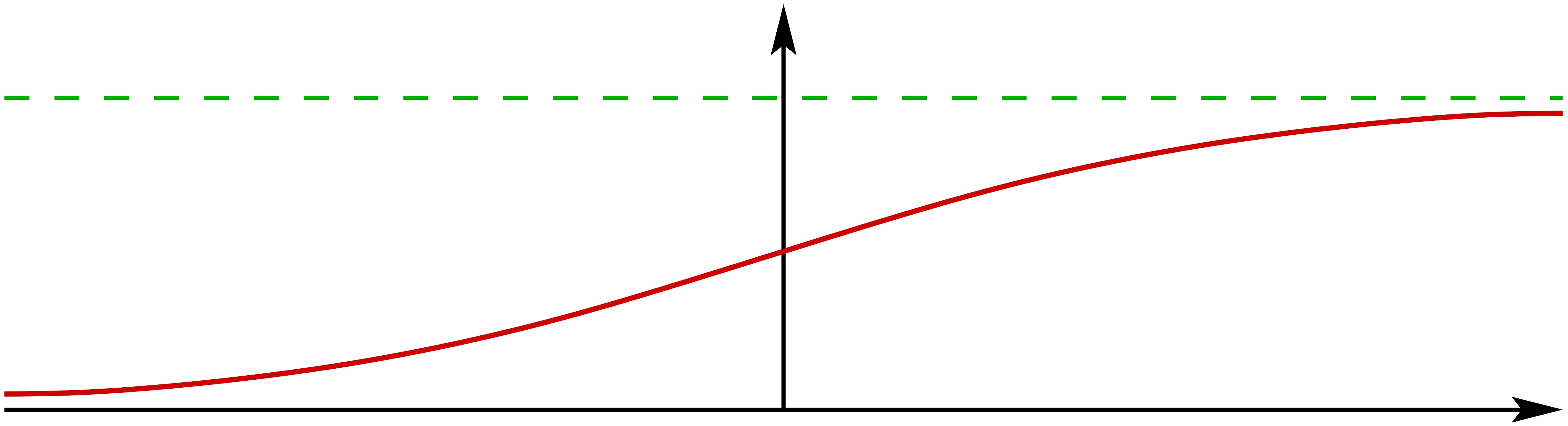}$\quad$
\includegraphics[width=0.4\textwidth]{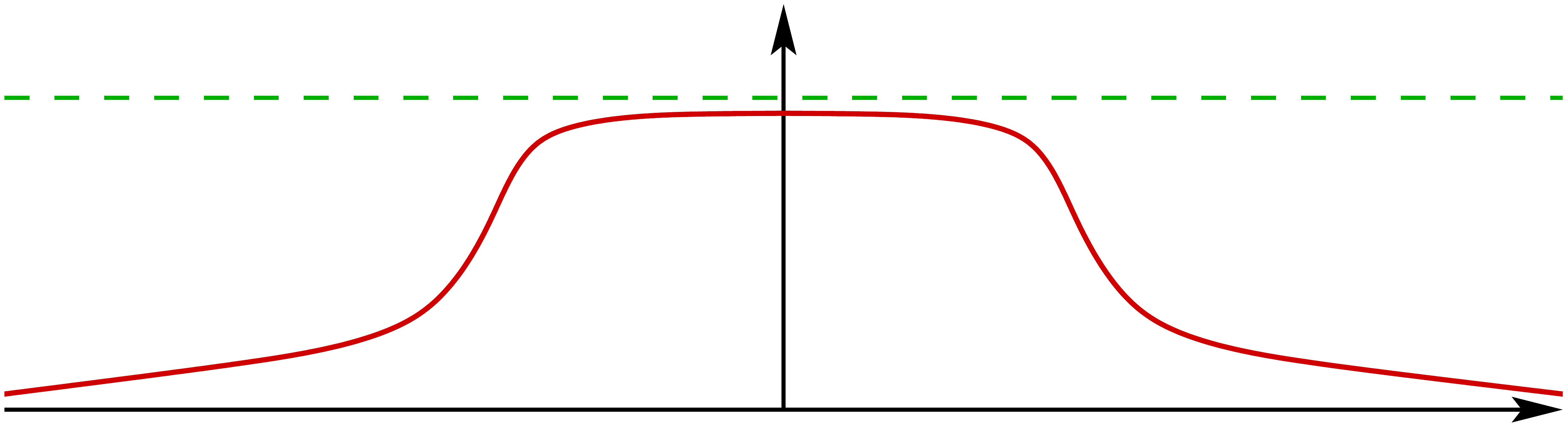}
\caption{A heteroclinic and a homoclinic connections.}\label{DEFG899}
\end{figure}

\begin{figure}[h]
\centering
\includegraphics[width=1\textwidth]{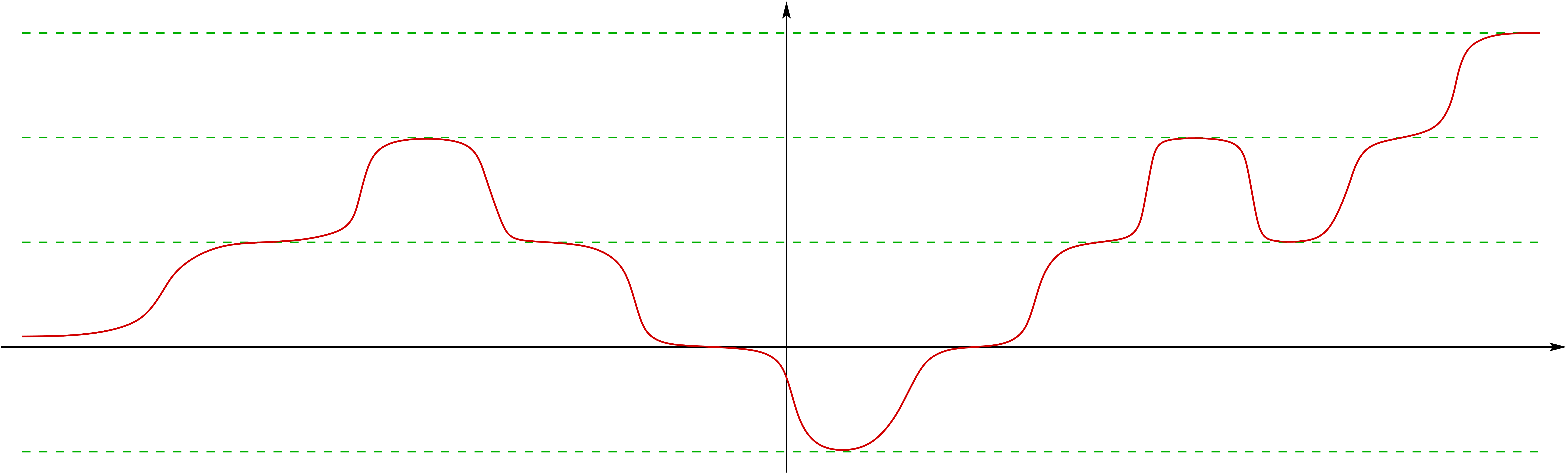}
\caption{A multibump orbit.}\label{DEFG8991}
\end{figure}

The goal in~\cite{MR3594365, MR4053239} is thus to
construct complex examples of stationary equilibria
for the atom dislocation
function induced by a perturbation of the potential
(roughly speaking, pressing or pulling a bit the crystal).
The complex example that we were interested in included:
\begin{itemize}
\item heteroclinic connections, i.e. transition layers connecting two minima of the potential~$W$,
\item homoclinic connections, i.e. a nontrivial connection of a minimum of the potential~$W$ to itself,
\item multibump orbits, i.e. solutions oscillating somewhat arbitrarily between
the minima of the potential~$W$,
\end{itemize}
see Figures~\ref{DEFG899} and~\ref{DEFG8991}.

To give a consistent mathematical framework to this perturbation setting, instead of~\eqref{SQRL3}
one considers the equation
\begin{equation}\label{SQRLPER}
(-\Delta )^s u(x)+a(x)\,W'(u(x))=0
\qquad{\mbox{for any }}\;\,x\in\R.\end{equation}
The role of~$a$ is that of a ``generic'' (arbitrarily small and with small derivatives,
but nontrivial) perturbation.

Among other results, it is proved in~\cite{MR3594365} that

\begin{theorem}\label{M12rSTHND}
Let~$s\in\left(\frac12,1\right)$, $N\in\N$, $\zeta_1,\dots,\zeta_N\in\Z$,
with~$|\zeta_{i+1}-\zeta_i|=1$,
for~$i\in\{1,\dots,N-1\}$.

Then, there exist~$b_1<\dots< b_{2N-2}\in\R$
and a solution~$u$ of~\eqref{SQRLPER}
such that:
\begin{eqnarray*}
&&\lim_{x\to-\infty} u(x)=\zeta_1, \\
&& \sup_{ x\in(-\infty,b_1] } |u(x)-\zeta_1|\le \frac1{100}, \\
&& \sup_{x\in[b_{2i},b_{2i+1}]}|u(x)-\zeta_{i+1}|\le\frac1{100}\quad{\mbox{ for all }}
i\in\{1,..., N - 2\},\\
&& \sup_{ x\in[b_{2N-2},+\infty) } |u(x)-\zeta_N|\le \frac1{100} \\
{\mbox{and }}&&\lim_{x\to+\infty} u(x)=\zeta_N.\end{eqnarray*}
\end{theorem}

In this spirit, Theorem~\ref{M12rSTHND} constructs a spatially oscillatory equilibrium
for the dislocation function as the one depicted in Figure~\ref{DEFG8991}:
we stress that the levels~$\zeta_1,\dots,\zeta_N\in\Z$
can be prescribed arbitrarily as consecutive integers and therefore
Theorem~\ref{M12rSTHND} states that a generic perturbation can produce dislocation patterns
that oscillate essentially as one wishes; we recall that, in view of the notation stated
on page~\pageref{NOTDISDS}, a positive oscillation
of the dislocation function would correspond to an atom dislocation ``right'' with respect
to its rest position, and a negative oscillation
of the dislocation function would correspond to an atom dislocation ``left'' with respect
to its rest position and accordingly the oscillations ``up and down'' detected in
Theorem~\ref{M12rSTHND} would correspond to atoms arbitrarily dislocated to the left and to the right
of the lattice in the setting of Section~\ref{INTR}.
Being able to dislocate atoms at will opens the possibility of associating to
the atomic structure a ``symbolic dynamics'' that describes a state
of the system in terms of the corresponding minimum of the potential~$W$.

Particular cases of Theorem~\ref{M12rSTHND} include the situations displayed
in Figure~\ref{DEFG899}, namely:
\begin{itemize}
\item when~$N:=2$ and~$\zeta_2:=\zeta_1+1$
(or~$N:=2$ and~$\zeta_2:=\zeta_1-1$), one obtains a heteroclinic orbit, \item
when~$N:=3$ and~$\zeta_3:=\zeta_1$, one obtains a homoclinic orbit.\end{itemize}
We observe that while heteroclinic orbits also exist in case of a constant
perturbation~$a$, the existence of homoclinic and multibump orbits heavily depends on the fact that~$a$
is supposed to be generic (and in particular nonconstant):
indeed, heteroclinic orbits for constant~$a$ have been constructed in~\cite{MR3081641, MR3280032}
but it has been shown in~\cite{MR4108219} that
no homoclinic connection and no multibump solution
are possible when~$a$ is constant.

The result in Theorem~\ref{M12rSTHND} is actually flexible enough to
deal also with more general interaction kernels
and with equilibria not necessarily disposed on a regular lattice.
The result also applies to systems (though in this case the subsequent
layer~$\zeta_{i+1}$ is not completely arbitrary given~$\zeta_i$, since it must be chosen
as a suitable action minimizers 
among close neighbors).

In terms of methodology, while Theorem~\ref{M12rSTHND} takes its
inspiration from classical variational methods, constrained minimization techniques and
energy estimates, its proof needs to account for some of the difficulties
provided by the nonlocal setting of equation~\eqref{SQRLPER}.
In particular:
\begin{itemize}
\item the ``cut and paste methods'' used in dynamical systems to create competitors for
the energy functional and obtain energy estimates
are affected in our case by nonlocality
(for instance, the energy of the union of two trajectories is not the sum of their energies),
\item the long tails of the solutions provide significant energy contributions
(suggestively,
as in what Einstein called
Mach's principle,
``mass out there influences inertia here'' and ``local physical laws are determined by the large-scale structure of the universe''~\cite{MR0424186}).\end{itemize}

\begin{figure}[h]
\centering
\includegraphics[width=1\textwidth]{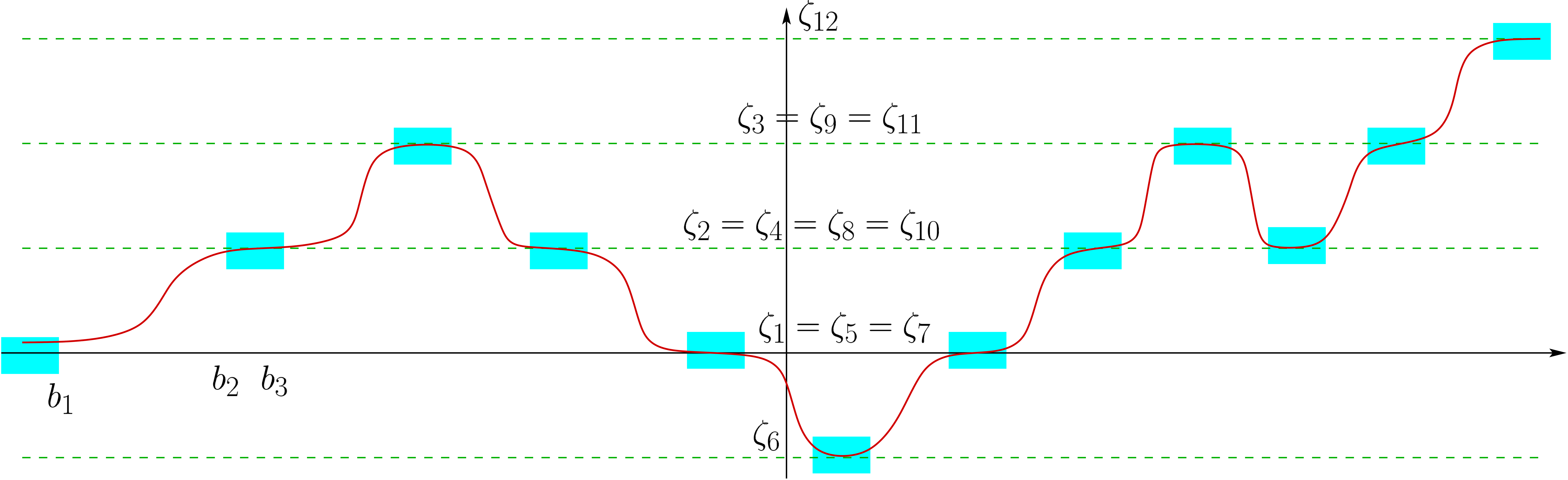}
\caption{Constrained minimization to construct a multibump orbit.}\label{MDEFG8991}
\end{figure}

To prove Theorem~\ref{M12rSTHND}, the variational method relies on a constrained minimization
problem: roughly speaking, one first minimizes the energy functional
under the additional request that the orbit passes through suitable ``windows''
strategically located to create the desired oscillations (namely, the cyan boxes in Figure~\ref{MDEFG8991}).

However, to make this argument work, one needs to avoid that the minimizing orbit constructed
in this way touches the boundary of the window (indeed, if it does not, then small perturbations
around the orbits are admissible and we obtain a local minimizer, hence a solution of the problem
with the desired oscillations).

To get rid of the possible touching between the minimal orbit and the windows, one
needs to carefully choose the intervals with endpoints~$b_1,\dots, b_{2N-2}$
and exploit the oscillations of~$a$ to place the transitions in a favorable energy range.

\begin{figure}[h]
\centering
\includegraphics[width=0.5\textwidth]{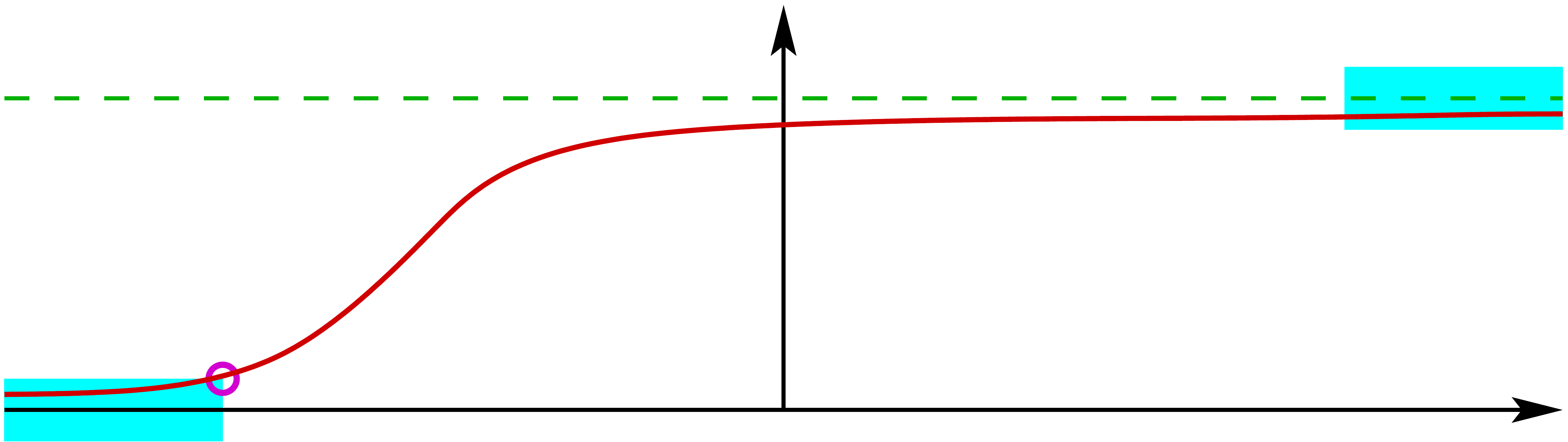}
\caption{The constrained minimization method going wrong.}\label{MDEFG8W991}
\end{figure}

To clarify this method, let us focus on the case of heteroclinics
(which will also serve as a cornerstone to build the other kinds of trajectories),
namely let us deal with the case~$N:=2$ and~$\zeta_2:=\zeta_1+1$.
Suppose that the method of constrained minimization ``goes wrong'', that is
a touching occurs as shown by the magenta circle in Figure~\ref{MDEFG8W991}.

One first shows that indeed the touching point must occur quite close to the boundary of the
window (at infinity, it is energetically more convenient for the solution
to stay close to the minima of the potential and thus away from the window).
Then, one needs to show that the transition from~$\zeta_1$ to~$\zeta_2$ occurs in a rather
short interval (it is energetically inconvenient for the solution to stay
away from the minima of the potential in a large interval).
Finally, once the solution gets sufficiently close to~$\zeta_2$ on the right, it is better
for it to remain close to it (it is energetically inconvenient for the solution to drift
away from the minima of the potential once it gets there).

\begin{figure}[h]
\centering
\includegraphics[width=0.5\textwidth]{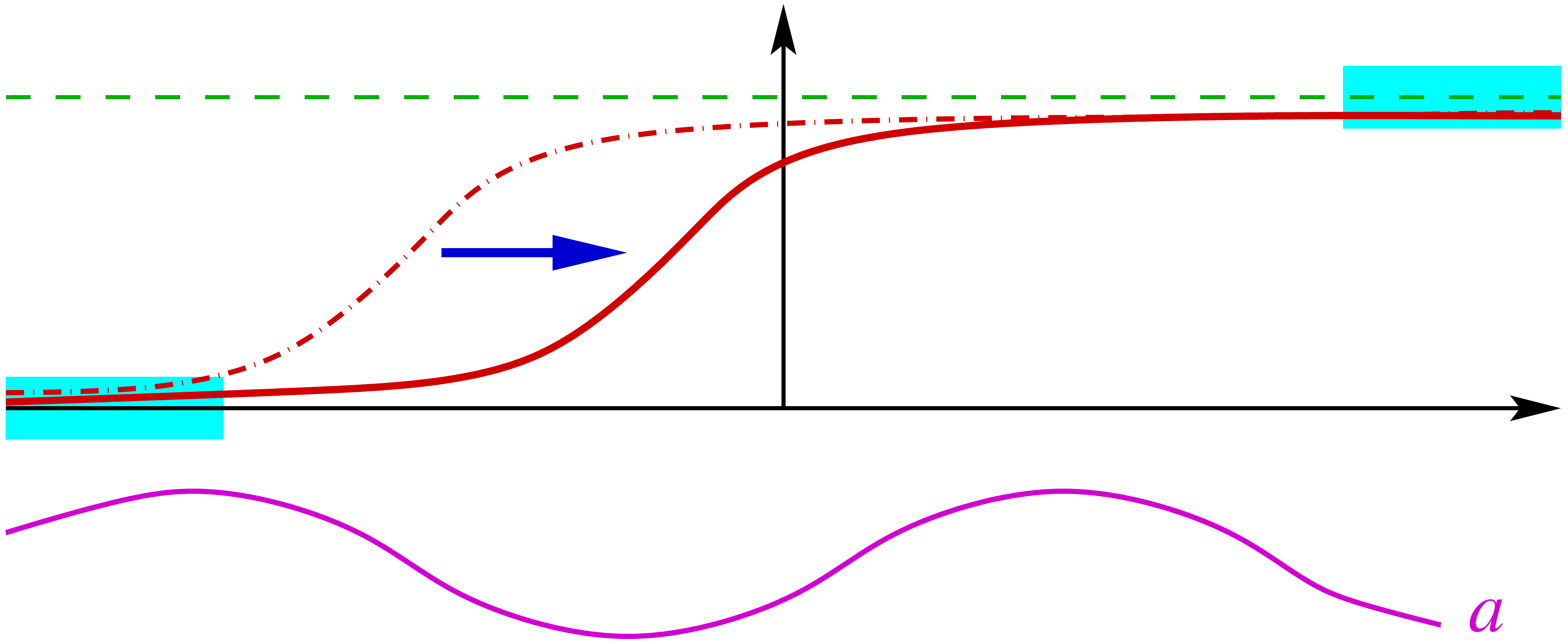}
\caption{Translating the orbit into an energetically more favorable position.}\label{MDEFG8W9912}
\end{figure}

With these three ingredients in mind, we realize that indeed
Figure~\ref{MDEFG8W991} represents essentially the ``worst case scenario''
that we want to avoid. For this, the idea is to translate the red trajectory in Figure~\ref{MDEFG8W991}
to the right and lower its energy: to accomplish this, one needs to assume that the modulating perturbation~$a$ presents some minima that we can enclose between the two cyan windows
(hence, to this end, if~$a$ is slowly oscillating one needs to place the windows sufficiently far from each other). Since the effect of~$a$ is to modulate the potential~$W$
and since we expect ``most of the energy'' to come from the transition from~$\zeta_1$ to~$\zeta_2$,
it will be energetically more convenient to locate the transition inside the minima of~$a$, thus
obtaining a new minimizer that does not touch the window, as desired (see Figure~\ref{MDEFG8W9912}).

When implementing these ideas, however, some effort is needed to take care of nonlocality,
since energy competitors cannot be obtained simply by gluing pieces of trajectories
(this feature will also become pivotal when one wants to address
Theorem~\ref{M12rSTHND} in its full generality, since homoclinic and multibump solutions
will be obtained precisely by ``almost'' glueing heteroclinics together).
To minimize the impact of nonlocality on the energy computations,
we glue orbits at ``very flat'' matching points, namely at points
at which two orbits meet with very small derivative, remaining close to each other
in a large interval (we quantify appropriately this notion and we
call such intervals ``clean intervals''). The existence of these clean intervals is
obtained by combining energy estimates and fractional
elliptic regularity theory.\medskip

\begin{figure}[h]
\centering
\includegraphics[width=0.3\textwidth]{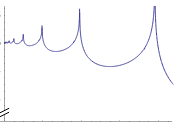}
\caption{A small energy competitor with many unbounded spikes.}\label{SPMDEFG8W9912}
\end{figure}

We observe that Theorem~\ref{M12rSTHND} was obtained in the fractional exponent range~$s\in\left(\frac12,1\right)$
which is arbitrarily close to, but different from, the exponent~$s=\frac12$ arising in~\eqref{SQRL1}.
The case~$s=\frac12$, or more generally~$s\in\left(0,\frac12\right]$,
is indeed quite challenging, due to several additional difficulties.
First of all, when~$s=\frac12$ one has to account for unbounded competitors, namely functions with finite energy
and unbounded spikes (which would make the constrained minimization problem with pointwise windows
unpractical). To appreciate this hindrance, one may consider the function
$$ \bar\psi(x):=\chi_{(-1,1)}(x)\ln(1-\ln|x|)$$
which belongs to~$ H^{1/2}(\R)\setminus L^\infty(\R)$, see~\cite[Appendix~A]{MR4053239}.
{F}rom this, one can easily construct competitors with small energy but unbounded spikes, e.g. both accumulating
at some given points and at infinity, such as
$$ \bar\psi_\e(x):=\e\sum_{k\in\Z} \frac{\bar\psi\big(2^{|k|}(x-2^k)\big)}{2^{|k|}},$$
see Figure~\ref{SPMDEFG8W9912}.

An additional difficulty comes from the fact that when~$s=\frac12$ heteroclinic solutions have infinite energy:
to appreciate\footnote{When~$s=\frac12$ and~$W(r)=\cos r$, \label{FOEX}
explicit heteroclinics are known, see e.g.~\cite[Appendices L and~M]{MR3967804}.
See also~\cite[Appendix~K]{MR3967804}
for more rigorous decay estimates than the ones sketched here.} this, one can argue that,
roughly speaking, if
\begin{equation}\label{LATYEBCIHBF}
{\mbox{$u_\star$ solves~\eqref{SQRL3} and is monotone and such that~$u_\star(-\infty)=0$ and~$u_\star(+\infty)=1$,}}
\end{equation} 
then, setting~$w:=u'_\star$, at~$+\infty$ one obtains a behavior of the type
$$ -(-\Delta )^s w(x)=W''(u_\star(x))w(x)\simeq W''(1) w(x),$$
that is, for large~$x$,
$$ w(x)\simeq\int_\R \frac{w(\zeta)-w(x)}{|x-\zeta|^{1+2s}}\,d\zeta
\simeq\int_\R \frac{w(\zeta)}{|x-\zeta|^{1+2s}}\,d\zeta\simeq\frac{C}{x^{1+2s}}
$$
and consequently
\begin{equation}\label{SLOWDDEC} u_\star(x)\simeq 1-\frac{C}{x^{2s}}.\end{equation}
This shows that when~$s=\frac12$, and more generally when~$s\in\left(0,\frac12\right]$, the heteroclinics possess\footnote{It
is instructive to remark that instead the potential energy is finite when~$s\in\left(\frac14,1\right)$ since in this range
\begin{eqnarray*}&& \int_\R W(u_\star(x))\,dx
=\int_{-\infty}^0\big( W(u_\star(x))-W(u_\star(-\infty)\big)\,dx+\int^{+\infty}_0\big( W(u_\star(x))-W(u_\star(+\infty)\big)\,dx\\&&\qquad
\le C\,\left(\int_{-\infty}^{-R}|u_\star(x)-u_\star(-\infty)|^2\,dx+\int^{+\infty}_R|u_\star(x)-u_\star(+\infty)|^2\,dx+1\right)\\&&\qquad
\le C\,\left(\int_{(-\infty,-R)\cup(R,+\infty)}\frac{dx}{|x|^{4s}}+1\right)<+\infty.
\end{eqnarray*}
}
infinite elastic energy, because, for large~$R$,
\begin{eqnarray*}&& \iint_{\R\times\R}\frac{|u_\star(x)-u_\star(y)|^2}{|x-y|^{1+2s}}\,dx\,dy
\ge \frac14\iint_{(R,+\infty)\times(-\infty,-R)}\frac{dx\,dy}{(x-y)^{1+2s}}
=\frac{1}{8s}\int_{(R,+\infty)}\frac{1}{(x+R)^{2s}}\,dx=+\infty.
\end{eqnarray*}
These observations clarify that for~$s=\frac12$, and more generally for~$s\in\left(0,\frac12\right]$,
the standard variational methods need to be modified since
the objects that we seek have infinite energy while the objects that we want to avoid (such as the ones presenting uncontrolled spikes)
can have finite, and even small, energy.

These difficulties have been addressed in~\cite{MR4053239},
where heteroclinic connections for equation~\eqref{SQRLPER} have been constructed
in the fractional exponent range~$s\in\left(\frac14,\frac12\right]$.
To overcome the difficulty arising from the infinite energy of the heteroclinics,
we use ``energy renormalization'' methods (by formally subtracting the infinite elastic energy of the heteroclinic
to the functional).
Also, to avoid finite energy competitors which may be discontinuous, and even unbounded,
we employe an ``energy penalization'' method, by
adding a ``viscosity'' term to the equation, obtaining uniform
estimates via elliptic theory (actually, via a combination of classical and fractional elliptic theories).
For this, an important step is to detect how the constrained minimizers 
detach from the windows:
namely, since the energy by itself cannot guarantee continuity in this fractional exponent range,
one needs an auxiliary argument to avoid that minimizers ``jump''  away when touching the boundary of the constraints.
These detachment estimates are obtained from the analysis of the corresponding obstacle problem
(the window playing the role of the obstacle at detaching points).

Another technical issue in this fractional exponent range is that
approximating heteroclinics may oscillate between equilibria
when placing the windows far from each other: to avoid this pathology, 
we use a second penalization method
to charge the $L^2$-difference with a centered excursion, so to control the ``baricenter''
of the heteroclinic.

\section{Dislocation dynamics}\label{DYN:SYS}

To understand the way in which dislocations move, one can look at the parabolic 
analogue of equation~\eqref{SQRL3}, namely consider solutions~$v=v(t,x)$ of
\begin{equation}\label{SQRL3-PARA}
\partial_t v=-(-\Delta )^s v-W'(v)\quad{\mbox{ in }}\,(0,+\infty)\times\R .
\end{equation}
More generally, one can consider the equation
$$ \partial_t v=-(-\Delta)^s v-W'(v)+\sigma_\e,$$
where~$\sigma_\e$ is a suitably small external stress, but to keep the exposition as simple as possible
we take here~$\sigma_\e\equiv0$.

In this section, we recall some results aiming at
describing the evolution of the parabolic dislocation function~$v$ for
suitably large space and time scales. To accomplish this goal, it is appropriate to
scale both the space and the time variables according to the homogeneity
of the equation, that is considering the function
\begin{equation}\label{SQRL3-PARA-eps00} v_\e(t,x):=v\left(\frac{t}{\e^{1+2s}},\frac{x}\e\right).\end{equation}
Here~$\e>0$ is a small parameter: sending~$\e\searrow0$ would thus
correspond to considering large time and space scales, in a suitably intertwined manner
that detects a coherent pattern. 

By scaling~\eqref{SQRL3-PARA}, one obtains the evolution equation for~$v_\e$ given by
\begin{equation}\label{SQRL3-PARA-eps}
\e\partial_t v_\e=-(-\Delta )^s v_\e-\frac1{\e^{2s}}W'(v_\e)\quad{\mbox{ in }}\,(0,+\infty)\times\R.
\end{equation}

Since the dynamics of the atom dislocations can be quite complicated,
it is also convenient to take into account ``well-prepared'' initial data,
obtained by the superposition of transition layers.
That is, we consider the heteroclinic solution~$u_\star$
constructed in~\cite{MR3081641, MR3280032} and mentioned in~\eqref{LATYEBCIHBF}
(normalized such that~$u_\star(0)=\frac12$): this will be the transition
layer (from~$0$ to~$1$) used as a building block for our construction. It is also convenient to consider a
``reversed'' transition layer (from~$1$ to~$0$) given by~$x\mapsto u_\star(-x)$.
Hence, we take as initial datum for~$v$ a superposition of transition layers and reversed
transition layers centered at given points~$\bar{x}_1<\dots<\bar{x}_N$.

{F}rom the physical point of view, one can imagine that the points~$\bar{x}_i$
represent the initial location of the atom dislocations:
the initial configuration is thus ``locally'' sufficiently close to an equilibrium, 
as described by the layer centered at~$\bar{x}_i$ (the other approaching a constant away
from their own center); notice that in this initial setting the atoms are located to either
the right or the left of their natural rest position, according to the choice of layer or reversed layer.

To formalize this aspect (that is, to distinguish between transition layers
and reversed transition layers), we let~$\zeta_1,\dots,\zeta_N\in\{-1,1\}$.
Roughly speaking, each~$\zeta_i$ denotes the ``orientation'' of the dislocation
located at~$\bar x_i$, with~$\zeta_i=1$ corresponding to
a layer transition and~$\zeta_i=-1$ to a reversed one.

With this notation, the initial datum for equation~\eqref{SQRL3-PARA-eps} is taken of the form
\begin{equation}\label{SQRL3-PARA-eps-6a}\bar v_\e(x):=\bar v_\e(0,x):=
\sum_{i=1}^N u_\star\left(\frac{\zeta_i}\e(x-\bar x_i)\right)-K,\end{equation}
where we denote by~$K\in\{1,\dots,N\}$ the number of negative orientations,
i.e.
$$ K:=\#\{i{\mbox{ s.t. }}\zeta_i=-1\}.$$
The reason for subtracting the constant~$K$ in~\eqref{SQRL3-PARA-eps-6a}
is that in this way if we have~$N$ even and~$K=N/2$ the function in~\eqref{SQRL3-PARA-eps-6a}
is normalized such that~$v_\e(-\infty)=v_\e(+\infty)=0$ (we will come back to this normalizing
choice on page~\pageref{SQRL3-PARA-eps-6aPA}).

The main result in this setting is that,
as~$\e\searrow0$,
the dislocations have the tendency to concentrate at single points of the crystal,
where the size of the slip coincides with the natural periodicity of the medium: more explicitly,
as~$\e\searrow0$,
we have that~$v_\e$ approaches a step function~$v_0$ of the form
\begin{equation}\label{SQRL3-PARA-eps-6} v_0(t, x):=\sum_{i=1}^N H\left(\zeta_i(x - x_i(t))\right)-K ,\end{equation}
where~$H$ is the Heaviside function (and we stress that the jump of the
Heaviside function is unitary, hence corresponding to the lattice periodicity).

The entity~$x_i(t)$ in~\eqref{SQRL3-PARA-eps-6} describes the discontinuities
of the ``limit dislocation function'' as the time~$t$ varies, that is, it outlines
the dynamics of the dislocation in time.

We have that the motion of these dislocation points is governed by an interior potential.
This potential is either repulsive or attractive,
depending on the mutual orientations
of the dislocations. More precisely, the dislocation dynamics described by~$x_i(t)$
correspond to the solutions of the Cauchy problem
of a system of ordinary differential equations of the form
\begin{equation}\label{SQRL3-PARA-eps-7}\begin{cases}
\dot x_i(t)=\displaystyle\gamma\,
\sum_{{1\le j\le N}\atop{j\ne i}}\zeta_i\zeta_j\frac{x_i(t)-x_j(t)}{2s\,|x_i-x_j|^{1+2s}}
,\\ x_i(0)=\bar{x}_i,
\end{cases}\end{equation}
where
\begin{equation}\label{3.6BIS}
\gamma:=\frac1{\| u_\star'\|_{L^2(\R)}}.\end{equation}

Interestingly, the system in~\eqref{SQRL3-PARA-eps-7}
is a gradient flow. We stress that when~$\zeta_i\zeta_j=1$, the potential
term is of repulsive type, but when~$\zeta_i\zeta_j=-1$ the term is attractive
and this may lead to collisions, as we will discuss in Section~\ref{COLL:SE}:
in particular, one can define~$T_c\in(0,+\infty]$ to be
the first time when a collision occurs.

Summarizing, by~\eqref{SQRL3-PARA-eps-6} 
and~\eqref{SQRL3-PARA-eps-7}, at an appropriately large space and time scale,
the behavior of the dislocations is governed by
an evolving step function (namely the one in~\eqref{SQRL3-PARA-eps-6})
and a system of ordinary differential equations describing its jump discontinuities
and accounting for the dynamics of the
dislocation points (as stated explicitly in~\eqref{SQRL3-PARA-eps-7}).

The precise mathematical formulation of this phenomenon goes as follows:

\begin{theorem}\label{PARTH}
For every~$\e > 0$, there exists a unique viscosity solution of~\eqref{SQRL3-PARA-eps}
in~$[0,+\infty)\times\R$,
with initial datum~\eqref{SQRL3-PARA-eps-6a}.

As~$\e\searrow 0$, this solution exhibits the following asymptotic behavior for all~$(t,x)\in
[0,T_c)\times\R$, for some~$T_c\in(0,+\infty]$:
\begin{equation}\label{KOTEPH2}
\begin{split}&
\limsup_{{(t',x')\to(t,x)}\atop{\e\searrow0}}v_\e(t',x')\le \overline{v_0}(t,x)\\
{\mbox{and }}\;&\liminf_{{(t',x')\to(t,x)}\atop{\e\searrow0}}v_\e(t',x')\ge \underline{v_0}(t,x),
\end{split}\end{equation}
where~$v_0$ is as in~\eqref{SQRL3-PARA-eps-6},
\begin{equation}\label{KOTEPH}
\begin{split}&
\overline{v_0}(t,x):=
\limsup_{{(t',x')\to(t,x)}}v_0(t',x')\\
{\mbox{and }}\;&
\underline{v_0}(t,x):=
\liminf_{{(t',x')\to(t,x)}}v_0(t',x')
.\end{split}\end{equation}
Moreover, the orbits~$x_i(t)$ in the definition of~$v_0$ in~\eqref{SQRL3-PARA-eps-6}
are solutions of the Cauchy problem in~\eqref{SQRL3-PARA-eps-7}.
\end{theorem}

We observe that outside the jumps~$x_i(t)$, the function~$v_0$ is continuous,
hence~\eqref{KOTEPH} boils down to
$$ \overline{v_0}(t,x)=\underline{v_0}(t,x)=v_0(t,x)\quad\mbox{ {provided that }}\;
x\not\in\{x_1(t),\dots,x_N(t)\},$$
and correspondingly~\eqref{KOTEPH2} reduces to
\begin{equation}\label{KOTEPH-3} \lim_{{(t',x')\to(t,x)}\atop{\e\searrow0}}v_\e(t',x')= {v_0}(t,x)
\quad\mbox{ {provided that }}\;
x\not\in\{x_1(t),\dots,x_N(t)\}.\end{equation}
In this sense, the statement in~\eqref{KOTEPH2} is a refinement
of the simpler one in~\eqref{KOTEPH-3} in which one takes into account the
upper and lower semi-continuous
envelopes of~$v_0$ to have a result comprising also the discontinuity set of the limit
dislocation function.

Theorem~\ref{PARTH} was first proved in~\cite{MR2851899} when~$s=\frac12$
(compare with footnote~\ref{FOEX}), and assuming that~$\zeta_1=\dots=\zeta_N=1$ (the case~$\zeta_1=\dots=\zeta_N=-1$ being analogous).

The case~$s\in\left(\frac12,1\right)$
with~$\zeta_1=\dots=\zeta_N=1$ has been established in~\cite{MR3296170}
and the case~$s\in\left(0,\frac12\right)$ with~$\zeta_1=\dots=\zeta_N=1$ 
in~\cite{MR3259559}.

Gradient flows as in~\eqref{SQRL3-PARA-eps-7} with~$\zeta_1=\dots=\zeta_N=1$ were studied in~\cite{MR2461827}.

Theorem~\ref{PARTH} in its full generality was then proved in~\cite{MR3338445},
where the case of collisions for the Cauchy problem in~\eqref{SQRL3-PARA-eps-7} under opposite orientations was also
taken into account.

\section{Collisions}\label{COLL:SE}

In light of the discussion in Section~\ref{DYN:SYS},
a natural question is to analyze the collision time~$T_c$ in Theorem~\ref{PARTH},
that is the first time in which two trajectories
of the system of ordinary differential equations in~\eqref{SQRL3-PARA-eps-7}
meet at the same point, namely
\begin{eqnarray*}&&
x_i(t)\ne x_j(t) {\mbox{ for all~$i\ne j$ when~$t\in[0,T_c)$,}}\\&&
x_{i_0}(T_c)=x_{i_0+1}(T_c)
{\mbox{ for some~$i_0\in\{1,\dots,N-1\}$.}}
\end{eqnarray*}
In particular, the system in~\eqref{SQRL3-PARA-eps-7} cannot be extended
for~$t>T_c$. We stress that~$T_c$ is finite only in presence of different orientations,
since when the orientations~$\zeta_i$ are all positive, or all negative,
then the potential is repulsive and particles do not collide.

The collision time~$T_c$ can be explicitly estimated in several concrete cases, such as the one
with two dislocations with opposite orientations, that of
three dislocations with alternate orientations, that
of~$N$ dislocations with alternate orientations, etc.

For instance,
among the several estimates obtained in~\cite{MR3338445}, we recall:

\begin{theorem}\label{COLLIS}
\noindent (i). Let~$N=2$ and~$K=1$. Then,
$$ T_c\le\frac{s(x\bar x_2-\bar x_1)}{(2s+1)\gamma},$$
where~$\gamma$ is as in~\eqref{3.6BIS}.
\smallskip

\noindent (ii). Let~$N = 3$, $\zeta_1=\zeta_2=1$ and~$\zeta_2=-1$. Let also
$$ \tau:=\frac{s\big(\min\{\bar x_2-\bar x_1,\,\bar x_3-\bar x_2\}\big)^{2s+1}}{(2s+1)\gamma}
\qquad{\mbox{and}}\qquad C:=\frac{2^{2s+1}}{2^{2s} - 1}
.$$
Then,
$$ T_c\in[\tau,\,C\tau].$$\smallskip

\noindent (iii) Let~$N\ge2$ and~$K\le N-1$, with~$\zeta_{i_0}\zeta_{i_0+1}=-1$ for
some~$i_0\in\{1,\dots,N-1\}$.
Then, there exists~$a_0\in\left(0,\frac1{(N-2)^{\frac1{2s}}}\right)$
such that the following statement holds true.

Assume that
$$ \bar x_{i_0+1}-\bar x_{i_0}\le a_0\min_{{i\in\{1,\dots,N-1\}}\atop{i\ne i_0}}\bar x_{i+1}-\bar x_{i}.$$
Then,
$$ x_{i_0+1}(t)-x_{i_0}(t)\le a_0\min_{{i\in\{1,\dots,N-1\}}\atop{i\ne i_0}}x_{i+1}(t)-x_{i}(t)$$
for all~$t\in[0,T_c)$, and
$$ T_c\le\frac{s(x\bar x_{i_0+1}-\bar x_{i_0})}{(2s+1)\gamma(1-(N-2)a_0^{2s})}.$$
\smallskip

\noindent (iv) Let~$N\ge2$ and~$\zeta_{i}\zeta_{i+1}=-1$ for
all~$i\in\{1,\dots,N-1\}$.

Then,
$$ T_c\le\begin{cases}\displaystyle
\frac{(N-1)(\bar x_N-\bar x_1)^{2s+1}}{(2s+1)\gamma} &{\mbox{ if $N$ is odd,}}\\ \\ \displaystyle
\frac{s(\bar x_N-\bar x_1)^{2s+1}}{(2s+1)\gamma} &{\mbox{ if $N$ is even.}}
\end{cases}$$
\end{theorem}

In particular, Theorem~\ref{COLLIS}(i) deals with the case of two dislocations with opposite
orientation and detects the collision time in terms of their mutual initial distance:
as a matter of fact, the result in~\cite{MR3338445} is more general than this, since
it takes into account also a possible external stress, in which case the result mentioned here
holds true provided that the dislocations are initially positioned sufficiently close to each other
(otherwise we also show that there are cases in which~$T_c=+\infty$).

Theorem~\ref{COLLIS}(ii) considers the case of three alternate dislocations
and provides bounds from above and below on the collision time. In this situation,
if two dislocations, say the first and the second, are closer than the others,
this property remains true for all times till collision occur:
namely, if~$\bar x_2-\bar x_1<\bar x_3-\bar x_2$ then~$x_2(t)-x_1(t)<x_3(t)-x_2(t)$
for all~$t\in[0,T_c)$.

This property is also useful to characterize ``triple collisions'', that is the situation
in which~$x_1(T_c)=x_2(T_c)=x_3(T_c)$. In the setting of Theorem~\ref{COLLIS}(ii),
one has that triple collisions occur if and only if the initial distance between the dislocation points is the same
and in this case one can characterize explicitly the collision time as the larger possible time
obtained in Theorem~\ref{COLLIS}(ii).
More precisely, triple collisions occur if and only if
$$ \bar x_2-\bar x_1=\bar x_3-\bar x_2
\qquad{\mbox{and}}\qquad T_c=C\tau.$$

Theorem~\ref{COLLIS}(iii) deals with the general case of~$N$ dislocations,
assuming that the orientation of the~$(i_0+1)$th and the one of the~$i_0$th are opposite.
In this situation, if the initial distance between these two dislocation points is sufficiently small
with respect to the others, then so it remains for all times till collision occurs,
and one can also give an explicit bound on the collision time.

Finally, Theorem~\ref{COLLIS}(iv) considers the case of~$N$ dislocations with alternate orientations
and bounds the collision time in terms of the maximal distance between initial dislocation points.

\section{Relaxation times and large time asymptotics}\label{sec:rel0686uyf}

Since the dawn of history, human beings have asked some fundamental questions, such as
``what happens after death?''.
Being unable to address this question, we focus on the weaker problem
of discussing what happens to the dislocations after collisions.

For this, we formally relate the collision time~$T_c$ discussed in Section~\ref{COLL:SE}
with the asymptotics presented in Theorem~\ref{PARTH}. Let us focus for simplicity
on the case of two dislocations with opposite orientations, i.e.~$N=2$ and~$K=1$,
say with~$\zeta_1=1$ and~$\zeta_2=-1$,
and let~$x_c:=x_1(T_c)=x_2(T_c)$ denote the collision point.
Thus, if~$x\ne x_c$ and~$t $ is sufficiently close to~$T_c$,
we have that~$x_1(t)\ne x$ and~$x_2(t)\ne x$,
whence it follows from~\eqref{SQRL3-PARA-eps-6} and~\eqref{KOTEPH-3} that,
when~$t$ is sufficiently close to~$T_c$,
\begin{equation*}
\lim_{{\e\searrow0}}v_\e(t,x)= {v_0}(t,x)
=\sum_{i=1}^N H\left(\zeta_i(x - x_i(t))\right)-1=
H(x - x_1(t))+H(x_2(t)-x)-1.\end{equation*}
In this setting, when~$t$ is sufficiently close to~$T_c$,
we have that either~$x<x_1(t)<x_2(t)$ 
(giving that~$H(x - x_1(t))=0$ and~$H(x_2(t)-x)=1$)
or~$x_1(t)<x_2(t)<x$
(giving that~$H(x - x_1(t))=1$ and~$H(x_2(t)-x)=0$).
{F}rom these observations we conclude that
\begin{equation*}
\lim_{{\e\searrow0}}v_\e(t,x)=0\end{equation*}
as long as~$t$ is sufficiently close to~$T_c$, whence
\begin{equation}\label{MEMSPD-1}
\lim_{t\nearrow T_c}\lim_{{\e\searrow0}}v_\e(t,x)=0.\end{equation}
However, in~\cite{MR3338445} it is also proved that
\begin{equation}\label{MEMSPD-2}\limsup_{{t\nearrow T_c}\atop{\e\searrow0}}v_\e(t,x_c)\ge1.\end{equation}
This says that while, by~\eqref{MEMSPD-1},
the two dislocations with opposite orientation average
after the collision, the dislocation function at the collision point and at the collision time
``keeps a bit of a memory'' of the collision event, as given by~\eqref{MEMSPD-2}.\medskip

In view also of this phenomenon,
it is interesting to
detect the ``relaxation time'' of the dislocation after collision. This question was investigated in~\cite{MR3511786}
in which it is established that, after a small transition time subsequent to the collision,
the dislocation function relaxes exponentially fast to a steady state. 
In this sense, the memory effect at the collision time detected in~\eqref{MEMSPD-2}
fades very fast.
Moreover, the system exhibits two separate scales of decay behaviors, namely an exponentially fast decay in time
and a polynomial decay in the space variable.
It is interesting that these two scales are somewhat asymptotically kept separate:
roughly speaking, for large times the main contribution comes from the time derivative of the dislocation
function balancing the potential (which produces an ordinary differential equation in time
leading to exponential decay), while in the space variable the leading term comes from the balance
of the fractional Laplacian of the dislocation function  with the potential (which, as in~\eqref{SLOWDDEC},
provides a slow decay of polynomial type).

\begin{figure}[h]
\centering
\includegraphics[height=3.8cm]{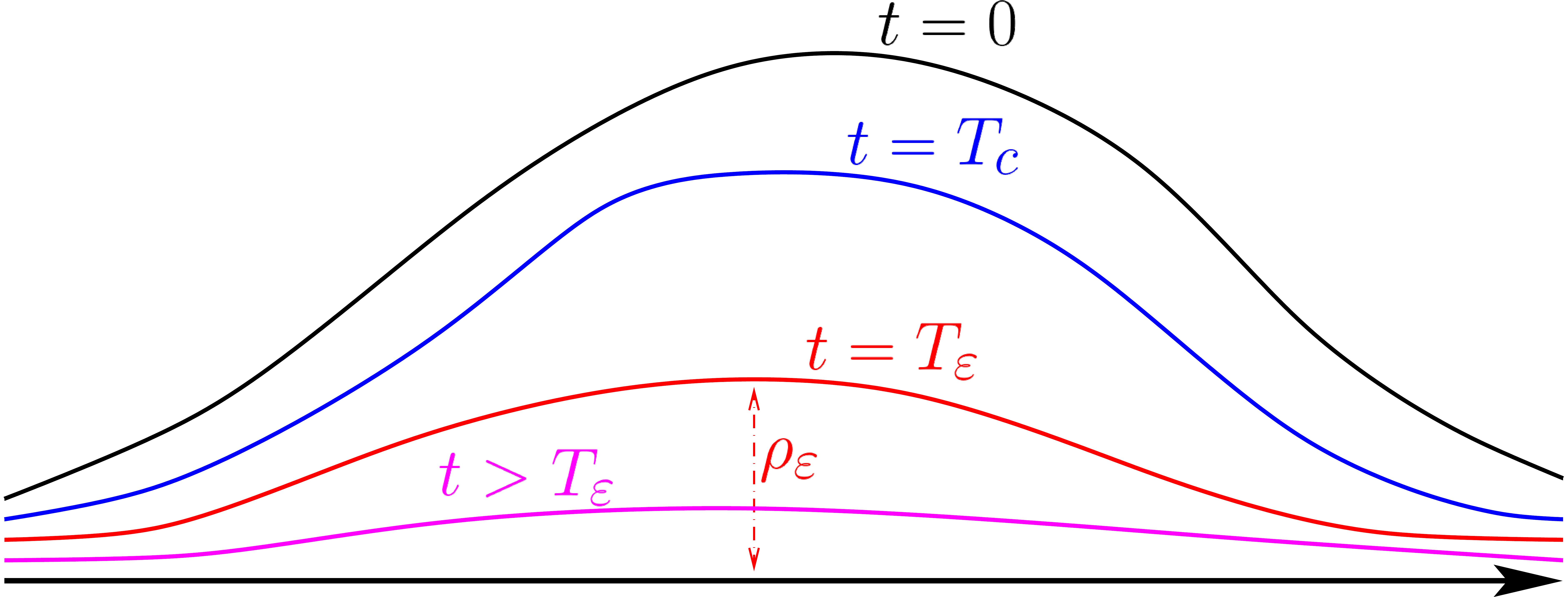}
\caption{The two dislocations situation described in Theorem~\ref{D2}.}\label{D2DEFE}
\end{figure}

The precise result for the two dislocations collision go as follows (see Figure~\ref{D2DEFE}
for a sketch of this situation):

\begin{theorem}\label{D2}
Let~$N=2$ and~$K=1$.
Then there exist~$\e_0> 0$ and~$c > 0$ such that for all~$\e\in(0,\e_0)$ we have that, for all~$t\ge T_\e$,
$$ |v_\e(t,x)|\le \rho_\e\exp\left(
-\frac{c}{\e^{2s+1}}(t-T_\e)
\right),$$
where~$\rho_\e\searrow0$ and~$T_\e\searrow T_c$ as~$\e\searrow0$.
\end{theorem}

\begin{figure}[h]
\centering
\includegraphics[height=3.8cm]{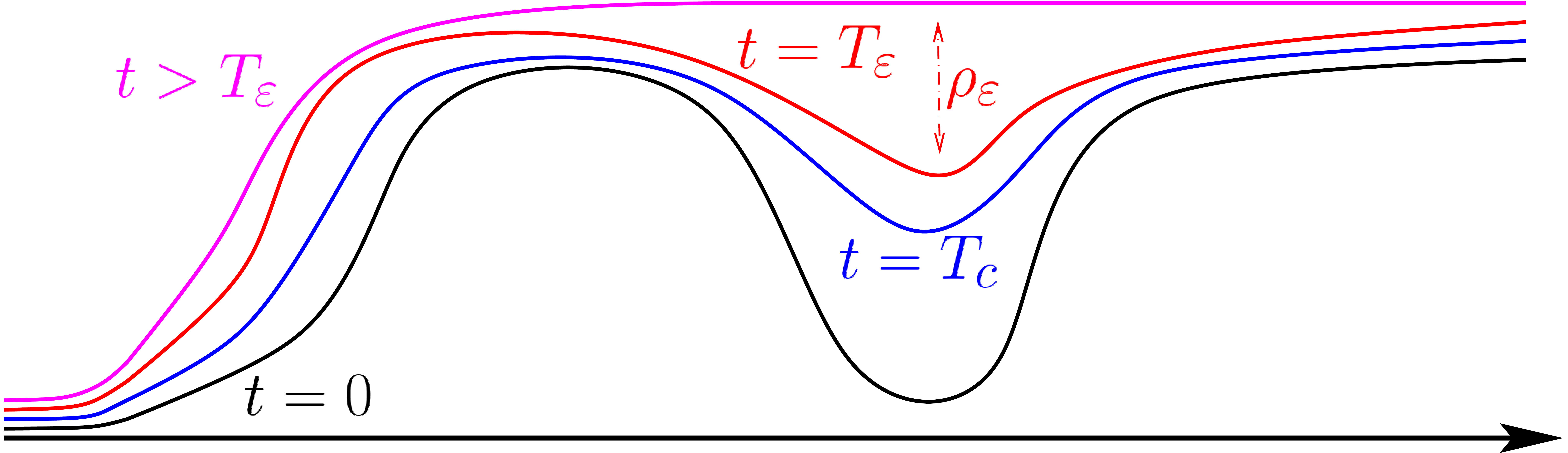}
\caption{The three alternate dislocations situation described in Theorem~\ref{D3}.}\label{D2DEFE4}
\end{figure}

The case
of three alternate dislocations is similar, albeit different, since the
steady state associated with this case is a heteroclinic and not the zero
function. That is, while in the case of two dislocations
the collision leads to an ``annihilation to the trivial equilibrium'',
in the case of three dislocations
we have that two dislocations ``annihilate each other'', but the third one ``survives''
and dictates the asymptotic behavior
(see Figure~\ref{D2DEFE4} for a sketch of this situation).
The precise result goes as follows
and provides a trapping from above and below for the solution~$v_\e$
in terms of the heteroclinic connection~$u_\star$:

\begin{theorem}\label{D3}
Let~$N=3$, $\zeta_1=\zeta_3=1$ and~$\zeta_2=-1$.
Then there exist~$\e_0> 0$ and~$c > 0$ such that for all~$\e\in(0,\e_0)$ we have that,
for all~$t\ge T_\e$,
$$ v_\e(t,x)\le 
u_\star\left(\frac1\e\left(x-y_\e+\kappa_\e\rho_\e
\left(1-\exp\left(
-\frac{c}{\e^{2s+1}}(t-T_\e)
\right)\right)\right)\right)+
\rho_\e\exp\left(
-\frac{c}{\e^{2s+1}}(t-T_\e)
\right)$$
and
$$ v_\e(t,x)\ge 
u_\star\left(\frac1\e\left(x-z_\e-\kappa_\e\rho_\e
\left(1-\exp\left(
-\frac{c}{\e^{2s+1}}(t-T_\e)
\right)\right)\right)\right)-
\rho_\e\exp\left(
-\frac{c}{\e^{2s+1}}(t-T_\e)
\right),$$
where~$\rho_\e\searrow0$, $\kappa_\e\searrow0$, $T_\e\searrow T_c$ and~$y_\e$, $z_\e\to x_\star$, for some~$x_\star\in\R$,
as~$\e\searrow0$.

In particular,
$$ \limsup_{t\to+\infty} v_\e(t,x)\le 
u_\star\left(\frac1\e\left(x-y_\e+\kappa_\e\rho_\e
\right)\right)$$
and
$$ \liminf_{t\to+\infty}v_\e(t,x)\ge 
u_\star\left(\frac1\e\left(x-z_\e-\kappa_\e\rho_\e
\right)\right)
.$$
\end{theorem}

The case of~$N$ dislocations has been discussed in~\cite{MR3703556}.
The general feature, besides technicalities, is that the parabolic character of the equation entails
a ``smoothing effect'' on the dislocation function slightly after collision, which leads, for large time,
to a precise asymptotics to either a constant function or a single heteroclinic,
depending on the algebraic properties of the orientations of the initial datum.

\begin{figure}[h]
\centering
\includegraphics[width=0.6\textwidth]{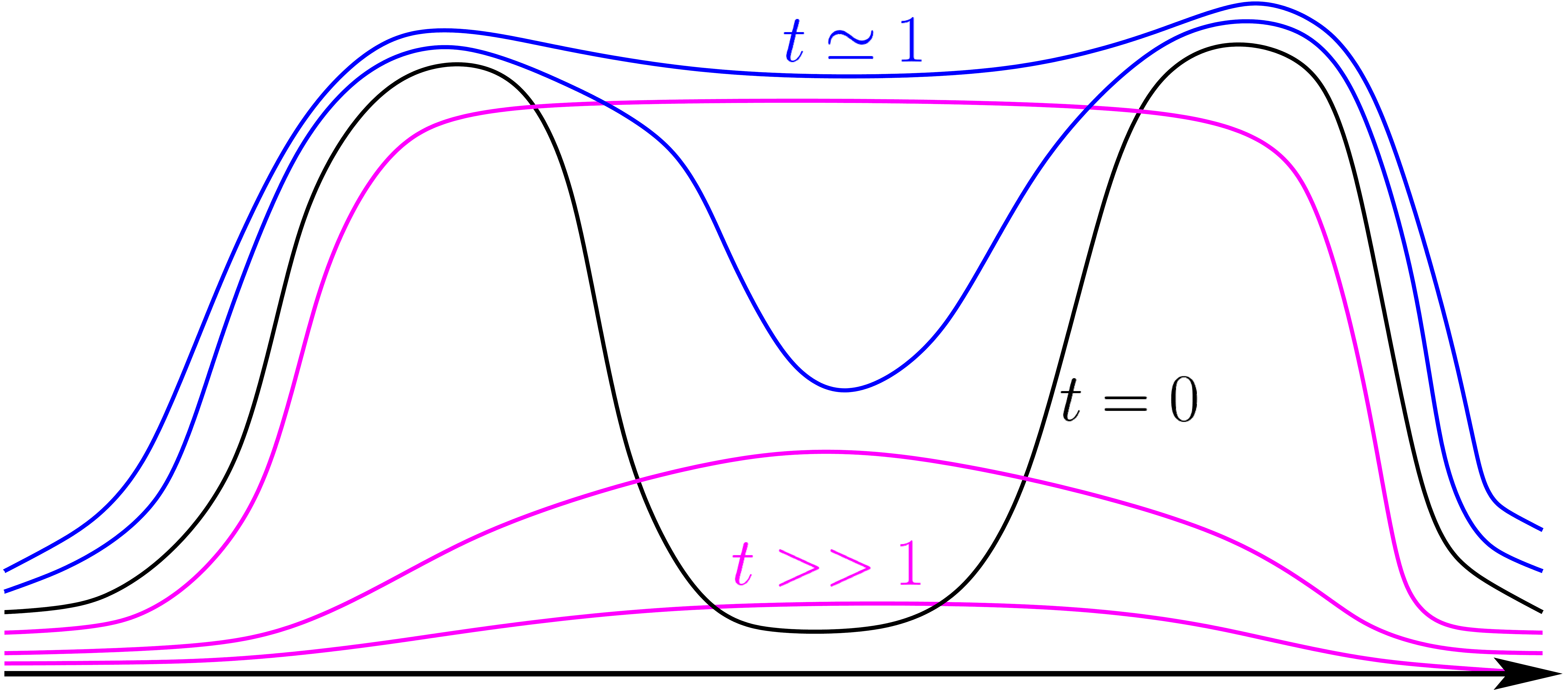}
\caption{Time asymptotic for the balanced case.}\label{D2DEFE-F1}
\end{figure}

To address the problem in its full complexity, it is convenient to distinguish
the case in which the dislocations are ``balanced'', i.e.~$N$ is even and~$K=N/2$,
corresponding to the situation in which the number of positively oriented dislocation
is equal to that of the negatively oriented ones,
and the case in which the dislocations are ``unbalanced'', namely
the number of positive ones exceeds that of negative ones (or vice versa).
A heuristic picture depicting the balanced case~$N=4$ and~$K=2$,
with~$\zeta_1=\zeta_3=1$ and~$\zeta_2=\zeta_4=-1$ is sketched in Figure~\ref{D2DEFE-F1}.

In the balanced regime, given the normalization in~\eqref{SQRL3-PARA-eps-6a}, \label{SQRL3-PARA-eps-6aPA}
we know that, the initial dislocation function~$\bar v_\e$ goes to zero both as~$t\to-\infty$
and as~$t\to+\infty$ and it turns out that
this asymptotic behavior in space also influences the asymptotic behavior in time.
More precisely, we have that after a transient time in which collisions occur, the dislocation function~$v_\e$
relaxes to zero exponentially fast.

The unbalanced case requires a more specific analysis. In this case, we can suppose
that~$K <N - K$ (the case~$K>N-K$ being similar) and define~$\ell:=N-2K>0$.
We observe that~$\ell$ counts the difference
between the positively oriented initial transitions (which are~$N-K$)
and the negatively oriented ones (which are~$K$). In this case, thanks to
the normalization in~\eqref{SQRL3-PARA-eps-6a}
the initial datum~$\bar v_\e$
is asymptotic to zero as~$x\to-\infty$ and to~$\ell$ as~$x\to+\infty$.

Given our previous discussion about the balanced case, one can expect that
the long time behavior of the dislocation function is influenced by these initial conditions
and one might believe that, as~$t\to+\infty$ the system will try to
average between the two values~$0$ and~$\ell$. However, this has to be taken with a pinch of salt,
since the constant~$\frac\ell2$ (that is the above mentioned average)
is not a minimum of the potential~$W$, unless~$\ell$ is even,
and therefore we cannot expect that~$v_\e$ converges to~$\frac\ell2$
as~$t\to+\infty$, unless~$\ell$ is even.

\begin{figure}[h]
\centering
\includegraphics[width=0.6\textwidth]{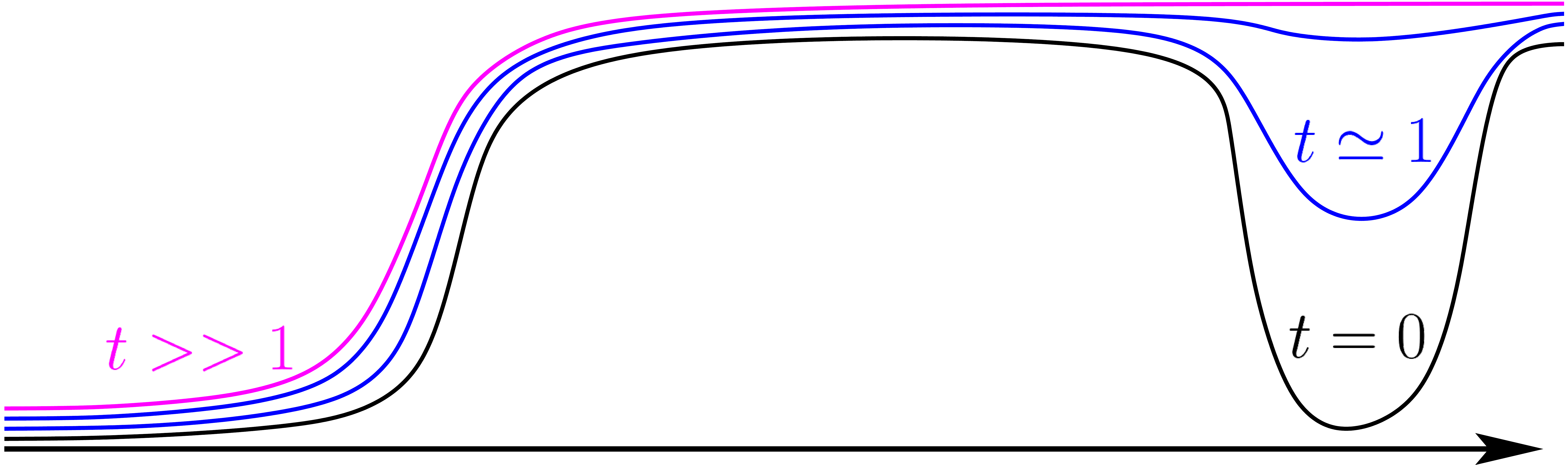}
\caption{Time asymptotic for the unbalanced case (with~$N=3$ and~$K=1$).}\label{D2DEFE-F2}
\end{figure}

The case~$\ell$ odd needs therefore an additional analysis.
On the one hand, in principle, the constant~$\frac\ell2$ could be a stationary point for~$W$,
hence it could provide a solution of the corresponding stationary equation in~\eqref{SQRL3}.
On the other hand, even if this happened, the
corresponding solution would be unstable from the
variational point of view, hence not suited to attract the dynamics of the parabolic flow.

As a result, when~$\ell$ is odd the average constant~$\frac\ell2$ is not the asymptotic limit of~$v_\e$
as~$t\to+\infty$, but it does serve as average asymptotics between the values~$\frac{\ell-1}2$
and~$\frac{\ell+1}2$, which are now integer and correspond to stable solutions of~\eqref{SQRL3}
induced by the minima of the potential~$W$.
What happens when~$\ell$ is odd is thus that, as~$t\to+\infty$, the dislocation function~$v_\e$
converges to a transition layer from~$\frac{\ell-1}2$
to~$\frac{\ell+1}2$.

For instance, a sketch of the situation in which~$N=3$ and~$K=1$, with~$\zeta_1=\zeta_3=1$
and~$\zeta_2=-1$ is given in Figure~\ref{D2DEFE-F2}. In this case, we have that~$\ell=3-2=1$
hence~$\ell$ is odd and the limit configuration in time is a transition from~$\frac{\ell-1}2=0$
to~$\frac{\ell+1}2=1$.

Instead, a sketch of the situation in which~$N=4$ and~$K=1$, with~$\zeta_1=\zeta_2=\zeta_3=1$
and~$\zeta_4=-1$ is given in Figure~\ref{D2DEFE-F2}. In this case, we have that~$\ell=4-2=2$,
which is even and the asymptotic state corresponds to the constant~$\frac\ell2=1$.
This asymptotics is indeed reached, roughly speaking, after the collision (and consequent 
relaxation) of the third and forth dislocation, followed by the repulsive interaction between
the first and second dislocation (this repulsion is responsible for the drift towards~$-\infty$
and~$+\infty$ of these two transition layers).

\begin{figure}[h]
\centering
\includegraphics[width=0.6\textwidth]{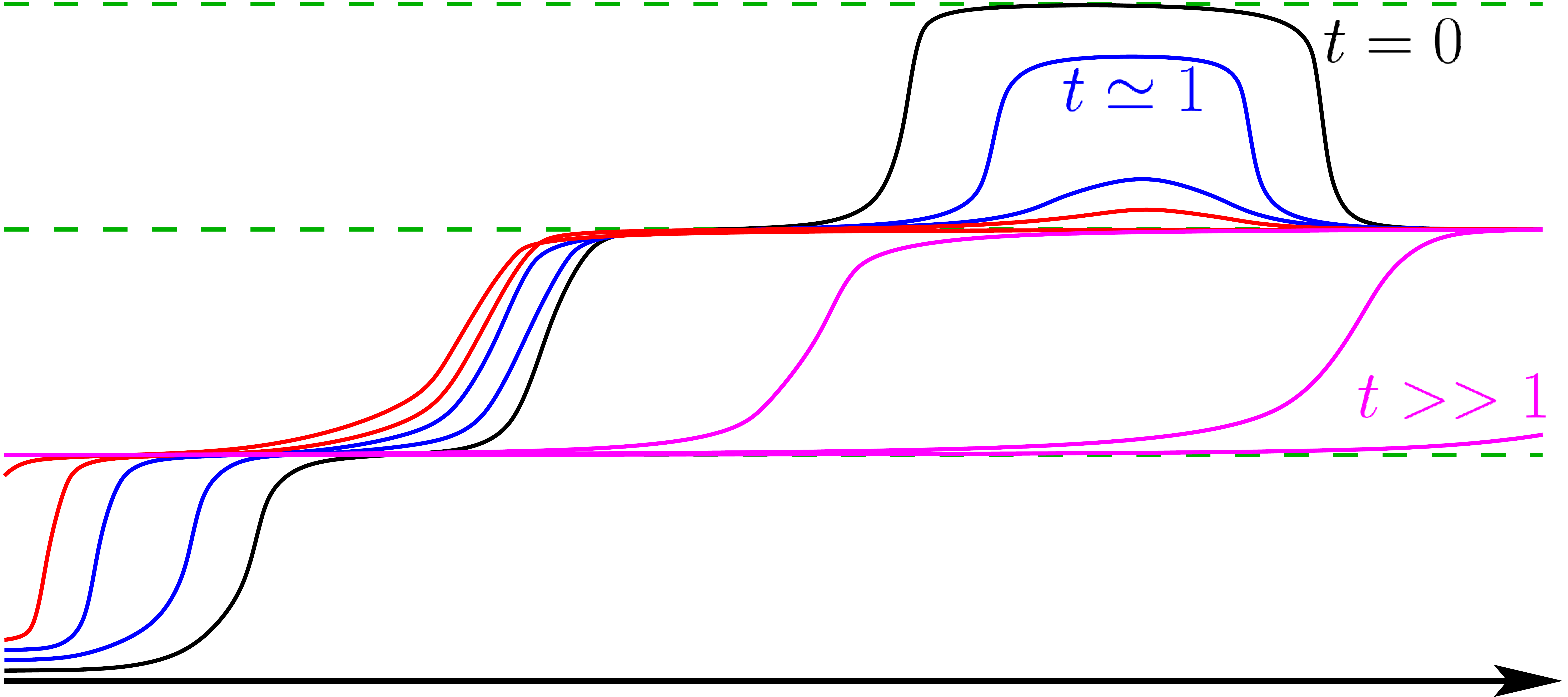}
\caption{Time asymptotic for the unbalanced case (with~$N=4$ and~$K=1$).}\label{D2DEFE-F3}
\end{figure}

See~\cite{MR3703556} for the corresponding formal statements.
We stress that the pictures of this paper are freehand drawn, not the outcome of a careful
numerical simulation, hence their intent is only explanatory,
they are at best qualitative and no attempt has been made to reproduce quantitatively
in the pictures
the different scales involved in the asymptotics.

\section{Orowan's Law}\label{orowan123}

In this section  we will focus on the physical case $s=\frac12$ and we will deduce the so-called
Orowan's Law, stating that the plastic strain velocity is proportional to the product of the density of dislocations by the effective stress.

In Section \ref{DYN:SYS} we have seen that the dynamics of $N$ parallel straight dislocation lines, all lying in the same slip plane, and all positive oriented ($\xi_i=1$ for $i=1,\ldots,N$) can be described by the 
ODE's system 
\begin{equation}\label{SQRL3-PARA-eps-7orowan}\begin{cases}
\dot x_i(t)=\displaystyle\gamma\,
\sum_{{1\le j\le N}\atop{j\ne i}}\frac{1}{x_i-x_j}
,\\ x_i(0)=\bar{x}_i,
\end{cases}\end{equation}
recall~\eqref{SQRL3-PARA-eps-7}
and~\eqref{3.6BIS}.

Here $T_c=+\infty$ as points with same orientation repeal each other, thus never collide. 
In this setting, we are  interested in the case when the number $N$ of dislocation points goes to infinity.
We will present below some results proven in  \cite{MR4138381} and aiming at identifying  at large scale an evolution model for the dynamics of a density of dislocations.
 For this, let  us introduce a new parameter $\delta=\delta_\e$ 
such that $\delta\to 0$ as $\e\to0$ and let us consider  the solution  $v_\e$ of 
\eqref{SQRL3-PARA-eps} with initial condition \eqref{SQRL3-PARA-eps-6a} (with $s=\frac12$ and $\xi_i=1$), where now $$N=N_\delta\sim \frac{1}{\delta}\to+\infty\quad\text{as }\e\to 0.$$
We next consider the following rescaling
$$u_\e(t,x):=\delta v_\e\left(\frac{t}{\delta},\frac{x}{\delta}\right).$$
Then $u_\e$ is solution to 
\begin{equation}\label{SQRL3-PARA-epsorowan}
\e\partial_t u_\e=-(-\Delta )^\frac{1}{2} u_\e-\frac1{\e}W'\left(\frac{u_\e}{\delta}\right)\quad{\mbox{ in }}\,(0,+\infty)\times\R
\end{equation}
with initial condition 
\begin{equation}\label{SQRL3-PARA-eps-6aorowan} u_\e(0,x)=
\sum_{i=1}^{N_\delta} \delta u_\star\left(\frac{x-\bar y_i}{\e\delta}\right)\quad \text{on }\R,\end{equation}
where, for $i=1,\ldots, N_\delta$, 
\begin{equation*}\bar y_i:=\delta \bar x_i\end{equation*}
and~$u_\star$ is the transition layer introduced in~\eqref{LATYEBCIHBF}.

Notice that with this rescaling $u_\e(0,x)$ (and thus, by comparison principle,
$u_\e(t,x)$) is bounded uniformly in $\e$ as 
$|u_\e(0,x)|\leq N_\delta \delta\sim1 $. 

The initial condition \eqref{SQRL3-PARA-eps-6aorowan} can be replaced, in further generality,  by any non-decreasing $C^{1,1}$ function $u_0$.
On the other hand,  any function $v\in C^{1,1}(\R)$ which is  non-decreasing and such that $v(-\infty)=0$, can be approximated by a function of the form 
\eqref{SQRL3-PARA-eps-6aorowan}.  
Indeed, given such a $v$,
for $0<\delta<1$,  define the $\delta i$ level set points $y_i$ as follows
\begin{equation}\label{yibar}y_i:=\inf\{x\in\R\,{\mbox{ s.t. }}\, v(x)=\delta i\}\quad {\mbox{ for all }}
i=1,\ldots, N_\delta,\end{equation}
where $$N_\delta:=\left \lfloor\frac{ v(+\infty)-\delta}{\delta} \right \rfloor.$$
Since the function~$v$ is non-decreasing, we have that $ y_i< y_{i+1}$. 
Then, $v$ can be approximated  in $L^\infty(\R)$ as stated in the following result.  
\begin{proposition}\label{approxpropfinal}
Let $v\in C^{1,1}(\R)$ be non-decreasing,  non-constant and such that~$v(-\infty)=0$.  Let 
$ y_i$ be defined as in \eqref{yibar}. 
Then, for all  $x\in\R$, 
$$
\left|\sum_{i=1}^{N_\delta}\delta  u_\star\left(\frac{x- y_i}{\e\delta}\right)-v(x)\right|\leq o_{\e,\delta}(1),
$$
where $o_{\e,\delta}$ is independent of $x$ and $o_{\e,\delta}(1)\to 0$ as $\e,\delta\to 0$.
\end{proposition}

Moreover, for $v$ and $y_i$ as above, the following approximation formula holds true: for any $i=1,\ldots,N_\delta$, 
\begin{equation}\label{i1onparticles}(-\Delta)^\frac{1}{2} v(y_{i})=\sum_{j\neq i}\frac{\delta}{y_i- y_{j}}+ o_\delta(1),\end{equation}
where $ o_\delta(1)\to 0$ as $\delta\to0$, uniformly on $\R$. 
Below is the main result in this setting.

\begin{theorem}\label{mainthmorowan}
Assume that $u_0 \in C^{1,1}(\R)$ and $u_0$ is  non-decreasing. Let $u_\e$ be the viscosity solution
of~\eqref{SQRL3-PARA-epsorowan}, with initial 
condition $u_\e(0,x)=u_0(x)$, for all~$x\in\R$.
Then,  as $\e\to0$, $u_\e$ converges locally uniformly in $(0,+\infty)\times\R$ to the unique  non-decreasing in $x$ viscosity solution~$\bar u$ of
\begin{equation}
\label{ubareq}\begin{cases}
\partial_t  u=-\gamma\partial_x  u(-\Delta)^\frac{1}{2 } u&\text{in }(0,+\infty)\times\R,\\
 u(0,\cdot)=u_0&\text{on }\R,
\end{cases}\end{equation}
where $\gamma$ is as in~\eqref{3.6BIS} 
\end{theorem}

The limit problem  \eqref{ubareq} represents the plastic flow rule for the   macroscopic crystal plasticity  with density of dislocations, where the solution
$\bar u$ can be interpreted as  the plastic strain,  $(-\Delta)^\frac{1}{2}\bar u$ is the internal stress created by the density of dislocations contained in a slip plane and identified by $\partial_x \bar u\ge 0$. 
Theorem~\ref{mainthmorowan} says that, in this regime, the plastic strain velocity 
is proportional to the dislocation density  times the effective stress. 
This physical law is  known as  Orowan's equation, see e.g. page~3739 in~\cite{sed}.

Equation 
\begin{equation}\label{limiteqorow}
\partial_t \bar u=-\gamma\partial_x \bar u(-\Delta)^\frac{1}{2}\bar u
\end{equation}
is an integrated form of a model studied by Head \cite{Head} for the self-dynamics of the dislocation
density~$\partial_x \bar u$.
Indeed, denoting by~$f:=\partial_x \bar u$, differentiating \eqref{limiteqorow}, we see that,  at least formally, $f$ solves
\begin{equation}\label{transporteq}\partial_t f=c_0\partial_x(f\mathcal{H}[f])\end{equation}
where $\mathcal{H} $is  the  Hilbert transform defined  by 
$$\mathcal{H}[f](x):=PV\int_\R\frac{f(y)}{y-x}dy.$$
Here we used that if $\bar u\in C^{1,\alpha}(\R)$ and $f=\partial_x \bar u\in L^p(\R)$ with $1<p<+\infty$, then   an integration by parts yields 
\begin{equation*}\label{halflaplhilbert}-(-\Delta)^\frac{1}{2}\bar u=\mathcal{H}[f].\end{equation*}
The conservation of mass  satisfied by  the positive integrable solutions of \eqref{transporteq} reflects the fact that if~$f=\partial_x \bar u$ is the density of dislocations, no dislocations are created or annihilated. 

Equation \eqref{transporteq}, called by Head  the  ``equation of motion of the dislocation continuum", was obtained by assuming that the speed
of each dislocation is proportional to the effective stress at the dislocation. 

This physical property is encoded in  the limit problem  \eqref{ubareq}. To see this, we first observe that 
the dislocation points can be approximated by the level set points of the limit function $\bar u$. More precisely,  define $y_i(t)$ as the $\delta i$-level set of $\bar u(t,\cdot)$ as in \eqref{yibar}.
Assuming  that $\bar u$ is smooth and $\partial_x \bar u (t, y_i(t))>0$, we see by differentiating in time  the equation $\bar u(t, y_i(t))=\delta i$ and using~\eqref{i1onparticles}
and~\eqref{ubareq}, that 
\begin{equation}\label{yidot}\dot y_i=-\frac{\partial_t \bar u(t, y_i(t))}{\partial_x \bar u(t, y_i(t))}=\gamma (-\Delta)^\frac{1}{2 } \bar u(t, y_i(t))\simeq \gamma \sum_{j\neq i}\frac{\delta}{ y_i- y_{j}}.\end{equation}
Scaling back to $x_i(t):=y_i(\delta t)/\delta$, we get 
$$\dot x_i \simeq \gamma \sum_{j\neq i}\frac{1}{ x_i- x_{j}},$$
which is, up to small errors,  \eqref{SQRL3-PARA-eps-7orowan}. Moreover, from \eqref{yidot}, we see  that the mean velocity of each  dislocation  is proportional to the effective stress $(-\Delta)^\frac{1}{2 } \bar u$ at the dislocation, as assumed by Head.

One can also take into account the case where dislocation points can be either positive or
negative oriented, by removing the assumption that the initial datum $u_0$ is non-decreasing.
In this situation Theorem \ref{mainthmorowan} is generalized in a forthcoming paper \cite{pnespre} as follows.

\begin{theorem}\label{mainthmcoll}
Assume that $u_0 \in C^{1,1}(\R)$ and that the limits $u(+\infty)$ and $u(-\infty)$ exist. Let $u_\e$ be the viscosity solution of~\eqref{SQRL3-PARA-epsorowan}, with initial 
condition $u_\e(0,x)=u_0(x)$, for all~$x\in\R$.
Then,  as $\e\to0$, $u_\e$ converges locally uniformly in $(0,+\infty)\times\R$ to the viscosity solution~$\bar u$
of
\begin{equation}
\label{ubareqcoll}\begin{cases}
\partial_t  u=-\gamma|\partial_x u|(-\Delta)^\frac{1}{2 } u&\text{in }(0,+\infty)\times\R,\\
 u(0,\cdot)=u_0&\text{on }\R,
\end{cases}\end{equation}
where $\gamma$ is as in~\eqref{3.6BIS}. 
\end{theorem}

Differentiating equation $\partial_t \bar u=-\gamma|\partial_x \bar u|(-\Delta)^\frac{1}{2 }\bar u$ yields the  two following equations for the positive and negative part of $f=\partial_x \bar u$, $$ \partial_t f^+=\partial_x( f^+\mathcal{H}(f^+-f^-)),\quad \partial_t f^-=-\partial_x( f^-\mathcal{H}(f^+-f^-))$$ 
which are  the 1-D version of the 2-D  Groma-Balogh equations \cite{groma},   
a  macroscopic model 
describing the evolution of the density of positive and negative oriented parallel straight dislocation lines. 
In this model in which collisions may occur, the  
 mass of the dislocation density $|\partial_x \bar u|=f^++f^-$ is not in general preserved, while it is preserved the mass of $f^+-f^-$, more precisely: for all $t>0$, 
$\bar u(t,+\infty)=u_0(+\infty)$ and $\bar u(t,-\infty)=u_0(-\infty)$, as shown in \cite{pnespre}. 

\section{Homogenization}\label{homog123}

In this section,
we consider equation \eqref{SQRL3-PARA-epsorowan} in any dimensions $n\geq 1$
and we deal with the corresponding homogenization problem.

To this end, we recall that
the parameter $\e$ in~\eqref{SQRL3-PARA-epsorowan}
describes the ratio between the microscopic  and the mesoscopic scale, where the discrete dynamics of dislocations is given by the ODE system~\eqref{SQRL3-PARA-eps-7}. Thus, recalling 
the initial condition \eqref{SQRL3-PARA-eps-6a}, the parameter~$\e$ can also
be interpreted as the distance between dislocation points at microscopic scale.

In this section, we freeze~$\e$, namely
we keep this distance fixed but, as in Section~\ref{orowan123},
we let the number of dislocation points go to infinity. Hence,
without loss of generality, we can assume $\e=1$.
Thus we deal with the problem 
\begin{equation}\label{SQRL3-PARA-ephom}
\begin{cases}
\partial_t u_\delta=-(-\Delta )^\frac{1}{2} u_\delta-W'\left(\displaystyle\frac{u_\delta
}{\delta}\right)&{\mbox{ in }}\,(0,+\infty)\times\R^n,\\
u_\delta(0,x)=u_0(x)&\text{ on }\R^n.
\end{cases}
\end{equation}
The case $n>1$, in particular $n=2$, allows us
to consider dislocation lines, all lying in the same slip plane, which can be curved and thus  for which the reduction of dimension argument described in the introduction cannot be applied. 

One could also eventually add on the right hand-side of the PDE in 
\eqref{SQRL3-PARA-ephom} a term of the form $\sigma(t/\delta,x/\delta)$, where $\sigma$ is periodic in both $t$ and $x$, representing 
 a stress created by the obstacles in the crystal or/and an applied exterior stress. For simplicity of presentation, as before we will take~$\sigma\equiv 0$.
 
The behavior of the solution $u_\delta$ of \eqref{SQRL3-PARA-ephom} has been investigated in \cite{monpatr1}. 
From a mathematical viewpoint, problem \eqref{SQRL3-PARA-ephom} is a homogenization problem, and
as usual in periodic homogenization (recall that the potential~$W$ is periodic), the limit equation is determined by a cell problem. In this case, such a problem is, for any $p\in\R^n$ and $L\in\R$, the following
\begin{equation}\label{w}
\begin{cases}
\partial_{\tau} w=-(-\Delta )^\frac{1}{2}w+L-W'(w+p\cdot y)&\text{in}\quad (0,+\infty)\times\R^n,\\
w(0,y)=0& \text{on}\quad
\R^n,
\end{cases}
\end{equation}
and one looks for some $\lam\in\R$ for which $w-\lam\tau$ is
bounded. The precise result goes as follows.

\begin{theorem}\label{ergodic}
For $L\in\R$ and $p\in\R^n$, there
exists a unique viscosity solution $w\in C((0,+\infty)\times\R^n)$ of~\eqref{w}
and there exists a unique $\lam\in\R$ such that $w$
satisfies 
$${\mbox{$\displaystyle\frac{w(\tau,y)}{\tau}$ converges towards $\lam$ as
$\tau\rightarrow+\infty$, locally uniformly in $y$.}}$$ The real
number  $\lam$ is denoted by $\overline{H}(p,L)$. The function
$\overline{H}(p,L)$ is continuous on $\R^n\times\R$ and
non-decreasing in $L$.
\end{theorem}

The function $\overline{H}(p,L)$ defined through the cell-problem \eqref{w}  is called effective Hamiltonian, and the 
homogenized limit problem is given by 
\begin{equation}\label{ueffett}
\begin{cases}
\partial_{t} u=\overline{H}(\nabla_x u,-(-\Delta )^\frac{1}{2}u)
&\text{in}\quad (0,+\infty)\times\R^n,\\
u(0,x)=u_0(x)& \text{on}\quad \R^n,
\end{cases}
\end{equation}
as stated below. 

\begin{theorem}\label{convergencehom}The solution
$u_\delta$ of \eqref{SQRL3-PARA-ephom} converges towards the unique viscosity solution $u^0$ of
\eqref{ueffett} locally uniformly in $(t,x)$, as~$\delta\searrow0$,
where $\overline{H}$
is defined in Theorem \ref{ergodic}.
\end{theorem}

As for the limit problems \eqref{ubareq} and \eqref{ubareqcoll}, the function $u^0$ can be interpreted  as a macroscopic plastic strain
satisfying the macroscopic  plastic flow rule (\ref{ueffett}). 
Moreover  $(-\Delta )^\frac{1}{2}u^0$ is the stress created by the macroscopic density of dislocations, $\nabla u^0$ is the vectorial dislocation density and 
 $|\nabla u^0|$ is the scalar dislocation density.

 The behavior of $\overline{H}(p,L)$ for small  stress $L$  and small density $|p|$, in dimension $n=1$, has been investigated in \cite{monpatr2}.
In this regime is it shown that
\begin{equation}\label{or}
\overline{H}(p,L) \simeq \gamma |p| L,
\end{equation}
thus recovering again the Orowan's Law discussed in Section~\ref{orowan123}. More precisely:

\begin{theorem}\label{hullprop}Let $p_0,\,L_0\in\R$. Then the function $\overline{H}$ defined in Theorem \ref{ergodic} satisfies
\begin{equation*}\frac{\overline{H}(\e p_0,\e L_0)}{\e^2}\rightarrow
\gamma |p_0|L_0\quad\text{as }\e\rightarrow0^+ ,\end{equation*}
where $\gamma$ is as in~\eqref{3.6BIS}. 
\end{theorem}

When the exterior force $\sigma$ does not vanish identically but has zero mean value, 
it is expected, but not proven,  a threshold phenomenon as in \cite{imr}, that is
$$\overline{H}(p,L)=0\quad \mbox{if}\quad |L|\quad \mbox{is small enough}.$$
This means more generally that this homogenization procedure
describes correctly the mechanical behaviour of the stress at
large scales, but keeps the memory of the microstructure in the
plastic law with possible threshold effects.

The results  contained in Theorems~\ref{convergencehom} and~\ref{hullprop} have been generalized in \cite{MR3334171} to the case where $(-\Delta)^\frac{1}{2}$ is replaced by 
$(-\Delta)^s$ with $s\in(0,1)$. 
The scaling of the system and the results obtained are different
according to the fractional parameter~$s\in(0,1)$.
Namely, when~$s>1/2$
the effective Hamiltonian ``localizes'' and it only depends
on a first order differential operator.
Instead, when~$s<1/2$, the non-local features are predominant
and the effective Hamiltonian involves
the fractional operator of order~$s$.
That is, roughly speaking, for any~$s\in(0,1)$,
the effective Hamiltonian is
an operator of order~$\min\{2s,1\}$, which reveals the stronger
non-local effects present in the case~$s<1/2$.

The  homogenization problems become 
\begin{equation}\label{ueps>1/2}
\begin{cases}
\partial_{t}
u_\delta=-\delta^{2s-1}(-\Delta)^su_\delta-W'\left(\frac{u_\delta}{\delta}\right)&\text{in}\quad (0,+\infty)\times\R^n,\\
u_\delta(0,x)=u_0(x)& \text{on}\quad \R^n,
\end{cases}
\end{equation}
for $s>\frac{1}{2}$, and
\begin{equation}\label{ueps<1/2}
\begin{cases}
\partial_{t}
u_\delta=-(-\Delta)^s u_\delta-W'\left(\frac{u_\delta}{\delta^{2s}}\right)+&\text{in}\quad (0,+\infty)\times\R^n,\\
u_\delta(0,x)=u_0(x)& \text{on}\quad \R^n,
\end{cases}
\end{equation}
for $s<\frac{1}{2}$.

Notice that the two different scalings formally coincide when $s=\frac{1}{2}$. 
The solution $u_\delta$ of \eqref{ueps>1/2} converges as $\delta\rightarrow 0$ to the solution $u^0$ of the homogenized problem 
\begin{equation}\label{ueffetts>1/2}
\begin{cases}
\partial_{t} u=\overline{H}_1(\nabla_x u)
&\text{in}\quad (0,+\infty)\times\R^n,\\
u(0,x)=u_0(x)& \text{on}\quad \R^n,
\end{cases}
\end{equation}
with an effective Hamiltonian $\overline{H}_1$ which does not depend on the fractional Laplacian
anymore, while the solution 
$u_\delta$ of \eqref{ueps<1/2} converges as $\delta\searrow 0$ to  $u^0$ solution of the following problem
\begin{equation}\label{ueffetts<1/2}
\begin{cases}
\partial_{t} u=\overline{H}_2(-(-\Delta)^su)
&\text{in}\quad (0,+\infty)\times\R^n,\\
u(0,x)=u_0(x)& \text{on}\quad \R^n,
\end{cases}\end{equation}
 with an effective Hamiltonian $\overline{H}_2$  not depending  on the gradient. 

The functions  $\overline{H}_1$ and $\overline{H}_2$ are determined by the following cell problem: for $p\in\R^n$ and $L\in\R$, 
\begin{equation}\label{whoms}
\begin{cases}
\partial_{\tau} w=-(-\Delta)^s w+L-W'(w+p\cdot y)
&\text{in}\quad (0,+\infty)\times\R^n,\\ w(0,y)=0& \text{on}\quad
\R^n,
\end{cases}
\end{equation} 
 that is $\overline{H}_1$ and $\overline{H}_2$ are determined by the unique $\lam$ for which  $w-\lam\tau$ is bounded (according to the cases $s>\frac{1}{2}$ and 
$s<\frac{1}{2}$, respectively). 
More precisely, as before  we have:
\begin{theorem}\label{ergodichoms} For $L\in\R$ and $p\in\R^n$, there
exists a unique viscosity solution $w\in C((0,+\infty)\times\R^n)$ of~\eqref{whoms} and there exists a unique $\lam\in\R$ such that $w$
satisfies: 
\begin{equation*}
{\mbox{$\displaystyle\frac{w(\tau,y)}{\tau}$ converges towards $\lam$ as
$\tau\rightarrow+\infty$, locally uniformly in $y$.}}\end{equation*} The real
number  $\lam$ is denoted by $\overline{H}(p,L)$. The function
$\overline{H}(p,L)$ is continuous on $\R^n\times\R$ and
non-decreasing in $L$.
\end{theorem}
Once  the cell problem is solved, the convergence results go as follows:

\begin{theorem}\label{convergence>1/2}The solution
$u_\delta$ of \eqref{ueps>1/2} converges towards the solution $u^0$ of
\eqref{ueffetts>1/2} locally uniformly in $(t,x)$ as~$\delta\searrow0$,
where $$\overline{H}_1(p):=\overline{H}(p,0)$$ and $\overline{H}(p,L)$
is defined in Theorem \ref{ergodichoms}.
\end{theorem}

\begin{theorem}\label{convergences<1/2} The solution
$u_\delta$ of \eqref{ueps<1/2} converges towards the solution $u^0$ of
\eqref{ueffetts<1/2} locally uniformly in $(t,x)$ as~$\delta\searrow0$,
where $$\overline{H}_2(L):=\overline{H}(0,L)$$ and $\overline{H}(p,L)$
is defined in Theorem \ref{ergodichoms}.
\end{theorem}

The following extension of the Orowan's Law is also proven in \cite{MR3334171}.
\begin{theorem} Let 
$p_0,\,L_0\in\R$ with $p_0\neq 0$. Then the function $\overline{H}$ defined in Theorem \ref{ergodichoms} satisfies
\begin{equation*}\frac{\overline{H}(\delta p_0,\e^{2s} L_0)}{\e^{1+2s}}\rightarrow
\gamma|p_0|L_0\quad\text{as }\e\rightarrow0^+ ,\end{equation*}
where $\gamma$ is as in~\eqref{3.6BIS}. 
\end{theorem}

\section*{Acknowledgements}
The first and third authors are members of INdAM/GNAMPA and AustMS.
The first author is supported by
the Australian Research Council DECRA DE180100957
``PDEs, free boundaries and applications''. 
The third author is supported by
the Australian Laureate Fellowship
FL190100081
``Minimal surfaces, free boundaries and partial differential equations''.

\end{document}